\documentclass{article}
\usepackage{amsmath,amssymb,amsthm}

\newcommand{\bs}{\backslash}
\newcommand{\longto}{\longrightarrow}

\newcommand{\C}{\mathbb{C}}
\newcommand{\R}{\mathbb{R}}

\newcommand{\Z}{\mathbb{Z}}

\newcommand{\Hom}{\mathrm{Hom}}

\newcommand{\s}{\sigma}
\newcommand{\pphi}{\varphi}

\newcommand{\veps}{\varepsilon}
\newcommand{\calC}{\mathcal{C}}
\newcommand{\piS}{\pi_1(S^2\bs\{s_1,\,\ldots\, ,s_l\})}
\newcommand{\Mod}{\mathcal{M}_{\mathcal{C}}}
\newcommand{\taum}{\tau^{-}}
\newcommand{\w}{\omega}
\newcommand{\lag}{\mathcal{L}}
\newcommand{\re}{\mathrm{Re~}}
\newcommand{\im}{\mathrm{Im~}}
\newcommand{\pconj}{\calC_1\times\cdots\times\calC_l}
\newcommand{\calD}{\mathcal{D}}
\newcommand{\calDt}{\widetilde{\mathcal{D}}}
\newcommand{\pdbl}{\calD_1\times\cdots\times\calD_l}
\newcommand{\pdblt}{\calDt_1\times\cdots\times\calDt_l}

\theoremstyle{definition}
\newtheorem{defi}{Definition}[section]

\newtheorem*{rk}{Remark}
\newtheorem*{ack}{Acknowledgements}

\theoremstyle{plain}
\newtheorem{thm_intro}{Theorem}

\newtheorem{thm}{Theorem}[section]
\newtheorem{prop}[thm]{Proposition}
\newtheorem{cor}[thm]{Corollary}
\newtheorem{lem}[thm]{Lemma}

\title{Representations of the fundamental group of an $l$-punctured
sphere generated by products of Lagrangian involutions} 
\author{Florent Schaffhauser}
\makeatletter
\def\printnotation{{%
\def\indexname{Index of notation}
\begin{theindex}
\@input{\jobname.ntn}
\end{theindex}
}}
\makeatother
\makeglossary

\begin{document}

\maketitle

\begin{abstract}
In this paper, we characterize unitary representations of $\pi:=\piS$ whose
generators $u_1, \,\ldots\,, u_l$ (lying in conjugacy classes fixed initially)
can be decomposed as products of two Lagrangian involutions
$u_j=\s_j\s_{j+1}$ with $\s_{l+1}=\s_1$ . Our main result is that such
representations are exactly the elements of the fixed-point set of an
anti-symplectic involution defined on the moduli space
$\Mod:=\Hom_{\calC}(\pi,U(n))/U(n)$ . Consequently, as this fixed-point set is
non-empty, it is a Lagrangian submanifold of $\Mod$. To prove this, we use
the quasi-Hamiltonian description of the symplectic structure of $\Mod$ and
give conditions on an involution defined on a quasi-Hamiltonian $U$-space
$(M, \w, \mu: M \to U)$ for it to induce an anti-symplectic involution on
the reduced space $M//U := \mu^{-1}(\{1\})/U$ .
\end{abstract}

\renewcommand{\thefootnote}{}
\footnotetext{AMS subject classification~: 53D20, 53D30}
\footnotetext{keywords~: momentum maps, moduli spaces, Lagrangian
submanifolds, anti-symplectic involutions, quasi-Hamiltonian}
\renewcommand{\thefootnote}{\arabic{footnotetext}}
\addtocounter{footnote}{-2}

%\footnotetext[1]{AMS subject classification~: 53D20, 53D30}
%\footnotetext[2]{keywords~: moduli spaces,
% momentum maps, Lagrangian submanifolds, anti-symplectic involutions, quasi-Hamiltonian spaces}

\section{Introduction}

The fundamental group $\pi:=\pi_1(S^2\bs\{s_1,\,\ldots\, ,s_l\})$
\glossary{$\pi:=\pi_1(S^2\bs\{s_1,\,\ldots\, ,s_l\})$~: the fundamental group of an
$l$-punctured $2$-sphere} of an
$l$-punctured $2$-sphere ($l\geq 1$) has finite presentation
$<g_1,g_2,\,\ldots\, ,g_l~|~g_1g_2{\ldots}g_l=1>$, where $g_j$ stands for the
homotopy class of a loop around $s_j$. Therefore, giving a unitary
representation of this surface group (i.e. a group morphism $\rho$ from $\pi$ to
$U(n)$) amounts to giving $l$ unitary
matrices $u_1,u_2,\,\ldots\, ,u_l$ satisfying the relation $u_1u_2{\ldots}u_l=1$
(we always identify $\C^n$ with $\R^{2n}$ and endomorphisms of $\C^n$ with their matrices in the
canonical basis). One may then want to study representations with prescribed conjugacy
classes of generators~: given $l$ conjugacy classes $\calC =
(\calC_j = \{u\exp(i\lambda_j)u^{-1}\,:\,u\in U(n)\})_{1\leq j\leq l}$, 
do there exist $l$ unitary matrices $u_1,u_2,\,\ldots\, ,u_l$ satisfying
$u_j\in\calC_j$ and
$u_1u_2{\ldots}u_l=1$. The answer to this problem was given by Agnihotri
and Woodward in \cite{AW}, by Belkale in \cite{Be} and by Biswas in
\cite{Bi2}~: they gave necessary and sufficient conditions on the
$\lambda_j\in\R^n$ for the above question to have a positive
answer (before that, the case of $SU(2)$ was discussed in \cite{JW}, in
\cite{gal}, in \cite{KM} and in \cite{Bi1}). In the following, we will always focus our interest on
representations with prescribed conjugacy classes of
generators and denote by $\Hom_{\calC}(\pi,U(n))$ the set of such
representations (i.e. group morphisms $\rho:\pi\to U(n)$ such that
$\rho(g_j)\in\calC_j$ for all $j$). Coming back to the relation $u_1{\ldots}u_l=1$, one may
notice that if we decompose each rotation $u_j \in U(n)$ as a product of two
orthogonal symmetries (which are no longer unitary transformations, since
they reverse orientation, see section \ref{laginv} for a precise definition
of these symmetries) in the following way $u_1=\s_1\s_2$, $u_2=\s_2\s_3, \,\ldots\, ,
u_l=\s_l\s_1$, then the relation $u_1\,\ldots\,u_l=1$ is automatically satisfied, since
orthogonal symmetries are elements of order $2$. The appropriate orthogonal
symmetries to consider turn out to be orthogonal symmetries with respect to a
Lagrangian subspace of $\C^n$, which are just real lines of $\C$ when $n=1$.
A unitary representation of $\piS$ whose generators
$u_1,\,\ldots\, ,u_l$ admit a decomposition $u_j=\s_j\s_{j+1}$, where
$\s_j$ is the orthogonal symmetry with respect to a Lagrangian subspace
$L_j$ of $\C^n$, and where $\s_{l+1}=\s_1$, will be called a
\emph{Lagrangian} representation. The
natural question to ask is then the following one~: when is a given representation a
Lagrangian one~? Further, two unitary representations of $\pi$ with
respective generators $(u_1, \,\ldots\, ,u_l)$ and $(u'_1, \,\ldots\, ,u'_l)$ being
equivalent if there exists a unitary map $\pphi\in U(n)$ such that
$u'_j=\pphi u_j \pphi^{-1}$ for all $j$, what can one say about the
set of Lagrangian representations in the moduli space $\Mod
:=\Hom_{\calC}(\pi,U(n))/U(n)$ of unitary representations~?\\
\indent In this paper, we address these two questions. First, we denote by $L_0$ the horizontal
Lagrangian $L_0 := \R^n \subset \C^n$ of $\C^n$ and we call a representation $\s_0$-Lagrangian
if it is Lagrangian with $L_1=L_0$. We will see in subsection \ref{beta} that a
given represemtation is Lagrangian if and only if it is equivalent to a
$\s_0$-Lagrangian one. We then obtain the following characterization of
$\s_0$-Lagrangian representations~:

\begin{thm_intro}{\label{thm_invo}}
Given $l\geq 1$ conjugacy classes $\calC_1,\calC_2, \,\ldots\, ,\calC_l \subset U(n)$
of unitary matrices such that there exist $(u_1,u_2, \,\ldots\, ,u_l)\in
\calC_1\times\cdots\times\calC_l$ satisfying $u_1u_2{\ldots}u_l=1$, the
representation of $\piS$ corresponding to such a $(u_1,u_2,\,\ldots\,,u_l)$ is
$\s_0$-Lagrangian if and only if $u_l=u_l^t$, $u_{l-1}=\overline{u}_l^{-1}
u_{l-1}^t \overline{u}_l$,
\ldots\,, and $u_1=\overline{u}_l^{-1} \overline{u}_{l-1}^{-1}
{\ldots}\overline{u}_2^{-1} u_1^t
\overline{u}_2{\ldots}\overline{u}_{l-1} \overline{u}_l$~.
\end{thm_intro}

Theorem \ref{thm_invo} will be proved in subsection \ref{beta} (theorem
\ref{s0-Lag}). Second, we recall that
the moduli space $\Mod$ of unitary representations of the surface
group $\pi$ with prescribed conjugacy classes of generators is a
symplectic manifold (actually a stratified symplectic space, see \cite{LS}, since we
have to take into account the singularities in the manifold
structure, see subsections 2.4 and 6.2 in \cite{J}). This symplectic structure, first investigated in \cite{AB}
and in \cite{G}, can be obtained in a variety of ways (see for instance
\cite{GHJW,AMM,AM,MW} and the references therein). For our purposes,
we will use the one
given by Alekseev, Malkin and Meinrenken in \cite{AMM} and think of our
moduli space as a symplectic quotient obtained by reduction of a
quasi-Hamiltonian manifold. We then have the following description of
the set of equivalence classes of Lagrangian representations of $\pi$~: 

\begin{thm_intro}{\label{thm_moduli}}
The set of equivalence classes of Lagrangian representations of
$\pi=\piS$ is a Lagrangian submanifold of the moduli space
$\Mod=\Hom_{\calC}(\pi,U(n))/U(n)$ of
unitary representations of $\pi$ (in particular, it is always non-empty).
\end{thm_intro}

Theorem \ref{thm_moduli} will be proved in subsection \ref{lag_submanifold}
(theorem \ref{lag_moduli}).
The fact that there always exist Lagrangian representations was first
proved in \cite{FW}, where the dimension of the submanifold of
(equivalence classes of) Lagrangian representations was shown to be
half the dimension of the moduli space. For a proof of the non-emptiness using
ideas from (quasi-) Hamiltonian geometry, we refer to \cite{mythesis} or to
the forthcoming paper \cite{S} (see also theorem \ref{convexity}). For
now, we will use the quasi-Hamiltonian description of the symplectic
structure of the moduli space to prove theorem \ref{thm_moduli}.\\
\indent The main intuition to tackle the aforementioned problems is
the use of momentum maps to solve questions of linear algebra (see
\cite{K}), which first seemed relevant for this problem after studying
the case $n=2$ (see \cite{FMS}), and which fits right into place
with the important idea of thinking of the space of equivalence classes of
representations (that is, the moduli space $\Mod$) as a symplectic
quotient. In this framework, the key idea to solve our problem is to obtain
the set of Lagrangian representations as the
fixed-point set of an involution $\beta$, which is
first used to give the explicit necessary and sufficient conditions
for a representation to be $\s_0$-Lagrangian appearing in theorem
\ref{thm_invo} and then turns out to induce an anti-symplectic involution on
the moduli space.\\
\indent After reviewing some background material on Lagrangian
involutions (that will later explain how the involution $\beta$ is
obtained), we shall proceed with recalling the notion of
quasi-Hamiltonian space introduced in \cite{AMM} and then use it to
obtain the symplectic structure of the moduli space $\Mod$ (we will
restrict ourselves to representations of $\piS$ and give an explicit
description of the symplectic $2$-form in the case $l=3$). Then we will
show how to obtain Lagrangian submanifolds of a quasi-Hamiltonian
symplectic quotient, in a way which theorem \ref{thm_moduli} will later provide a concrete
example of. Finally, we will obtain $\s_0$-Lagrangian representations of
$\pi$ as the fixed-point set of an involution on the product
$\calC_1\times\cdots\times\calC_l$ of the prescribed conjugacy classes, and
therefrom deduce theorem \ref{thm_invo} and theorem \ref{thm_moduli}. Along
the way, we will have proved another result which is worth mentioning here
and which will be proved in section \ref{Lag_QH} (theorem \ref{lag_locus})~:

\begin{thm_intro}
Let $U$ be a compact connected Lie group and let $(M,\w)$ be a
quasi-Hamiltonian $U$-space with momentum map $\mu: M \to U$ . Let
$\tau$ be an involutive automorphism of $U$, denote by $\taum$ the involution
defined on $U$ by $\taum(u)=\tau(u^{-1})$ and let $\beta$ be an involution on
$M$ such that~: 
\begin{enumerate}
\item[(i)] $\forall u \in U$, $\forall x \in M$,
$\beta(u.x)=\tau(u).\beta(x)$
\item[(ii)] $\forall x \in M$, $\mu \circ \beta (x)=\taum \circ \mu (x)$
\item[(iii)] $\beta^* \w = -\w$
\end{enumerate}
then $\beta$ induces an anti-symplectic involution $\hat{\beta}$ on the
reduced space $M^{red}:=\mu^{-1}(\{1\})/U$ . If $\hat{\beta}$ has fixed
points, then $Fix(\hat{\beta})$ is a Lagrangian submanifold of $M^{red}$.
\end{thm_intro}

\begin{ack}
Before starting, I would like to thank my adviser Elisha Falbel for submitting the
above problem to me. Numerous discussions with him and with Richard
ãWentworth were of valuable help to me. I would also like to thank Alan
Weinstein for encouragement on the momentum map approach and Johannes
Huebschmann for mentioning the notion of quasi-Hamiltonian space to me.
My deepest gratitude goes to
Jiang Hua Lu and Sam Evens for that incredibly fruitful discussion we
have had in Paris in the Spring of 2004. It was a sincere pleasure.
I am also grateful to Professor Yoshiaki Maeda and the
department of Mathematics at Keio University in Yokohama for their
hospitality at the time this paper was being written. My presence in
Keio was made possible thanks to a short-term doctoral fellowship
granted by the Japan Society for the Promotion of Science (JSPS). Finally, I
would like to thank the referee for his comments and suggestions to improve
the readibility of this paper and for pointing out to me the results in
\cite{gal} and \cite{KM}.
\end{ack}

\section{Background on Lagrangian involutions and angles between Lagrangian
subspaces}\label{laginv}

We give here the properties of Lagrangian involutions that we shall need in
the following. Recall that $\C^n$ is endowed with the symplectic
form $\w=-\im h$ where $h$ is the canonical Hermitian product
$h=\sum_{k=1}^n dz_k \otimes d\overline{z}_k$, for which it is symplectomorphic
to $\R^{2n}$ endowed with the canonical symplectic form $\w=\sum_{k=1}^n
dx_k \wedge dy_k$ . Mutiplication by $i\in\C$ in $\C^n$ corresponds to an
$\R$-endomorphism $J$ of $\R^{2n}$ satisfying $J^2=-Id$. Denoting $g=\re
h=\sum_{k=1}^n (dx_k \otimes dx_k + dy_k \otimes dy_k)$ the canonical
Euclidean product on $\R^{2n}$, we have $g=\w(.,J.)$ ($J$ is called a complex
structure and is said to be compatible with $\w$). A real subspace $L$ of
$\C^n$ is said to be Lagrangian if $\w|_{L \times L}=0$ and if $\dim_{\R} L
= n$ (that is, $L$ is maximal isotropic with respect to $\w$). One may then
check that $L$ is Lagrangian if and only if its $g$-orthognal complement is
$L^{\perp_g}=JL$ . We may then define, for any Lagrangian subspace $L$ of
$\C^n$, the $\R$-linear map

$$\begin{array}{rrcl}
\s_L: & \C^n=L\oplus JL & \longto     & \C^n \\
      &      x + Jy     & \longmapsto & x-Jy
\end{array}$$

\noindent called the \emph{Lagrangian involution associated to} $L$. Observe that
$\s_L$
is \emph{anti-holomorphic}~: $\s_L \circ J = -J \circ \s_L$. In the
following, we denote by $\lag(n)$ \glossary{$\lag(n)$~: the Lagrangian
Grassmannian of $\C^n$} \glossary{$\s_L$~: the Lagrangian involution associated to
$L\in\lag(n)$} the set of all Lagrangian subspaces of
$\C^n$ (the Lagrangian Grassmannian of $\C^n$). Finally, recall that, under the
identification $(\C^n,h) \simeq (\R^{2n},J,\w)$, we have $U(n)=O(2n) \cap
Sp(n)$. Furthermore, the action of $U(n)$ on $\lag(n)$ is transitive and the
stabilizer of the horizontal Lagrangian $L_0:=\R^n\subset\C^n$ is the
orthogonal group $O(n) \subset U(n)$, giving the usual homogeneous
description $\lag(n)=U(n)/O(n)$ . Observe that $O(n)=Fix(\tau)$ where
$\tau : u \mapsto \overline{u}$ is complex conjugation on $U(n)$, so that $\lag(n)$ is a
compact symmetric space.

\begin{prop}\emph{\cite{FMS}}
Let $L\in\lag(n)$ be a Lagrangian subspace of $\C^n$. Then:
\begin{enumerate}
\item[(i)] There exists a unique anti-holomorphic map $\s_L$ whose
fixed point set is exactly $L$.
\item[(ii)] If $L'$ is a Lagrangian subspace such that
$\s_L=\s_{L'}$, then $L=L'$~: there is a one-to-one correspondence between
Lagrangian subspaces and Lagrangian involutions.
\item[(iii)] $\s_L$ is \emph{anti-unitary}~: for all $z,z'\in\C^n$,
$h(\s_L(z),\s_L(z'))=\overline{h(z,z')}$.
\item[(iv)] For any $\pphi \in U(n)$,
$\s_{\pphi(L)}=\pphi\s_L\pphi^{-1}$.
\end{enumerate}
\end{prop}

Denote then by $\mathcal{LI}nv(n):=\{\s_L~:~L\in\lag(n)\}$ the subset of $O(2n)$
consisting of Lagrangian involutions. Observe that it is not a subgroup, as
it is not stable by composition of maps. Statement (iv) of the above
proposition then shows that the subgroup $\widehat{U(n)}:=<U(n) \cup
\mathcal{LI}nv(n)>\,\subset O(2n)$ generated by Lagrangian involutions and
unitary transformations is in fact generated by $U(n)$ and $\s_{L_0}$~:
$\widehat{U(n)}=<U(n)\cup\{\s_{L_0}\}>$. As a word in
$<U(n)\cup\{\s_{L_0}\}>$ contains either an even or an odd number of
occurrences of $\s_{L_0}$ (depending only on whether it represents a
holomorphic or an
anti-holomorphic transformation of $(\R^{2n},J)\simeq\C^n$), it can be written
uniquely under the reduced form $u\veps$ where $u\in U(n)$ and $\veps=1$
or $\veps=\s_{L_0}$. Consequently, we have $<U(n)\cup\{\s_{L_0}\}>\,=U(n)
\sqcup U(n)\s_{L_0}$, so that $U(n)$ is indeed a subgroup of index $2$ of
$\widehat{U(n)}$. Further, if we write $\Z / 2\Z = \{1,\s_{L_0}\}$ and
consider the action of this group on $U(n)$ given by
$\s_{L_0}.u=\s_{L_0}u\s_{L_0}=\overline{u}=\tau(u)$, then the map
$$\begin{array}{rcl}
U(n)\rtimes \Z/2\Z & \longto & U(n) \sqcup U(n)\s_{L_0} \\
(u,\veps) & \longmapsto & u\veps
\end{array}$$ \noindent (where $\veps=1$ or $\veps=\s_{L_0}$) is a group
isomorphism. Finite subgroups of $U(2)\rtimes\Z/2\Z$ generated by Lagrangian
involutions are studied in \cite{F}. As for us, one of the major interests of
Lagrangian involutions will be that they measure angles of Lagrangian subspaces
of $\C^n$ under the action of the unitary group~:

\begin{thm}\emph{\cite{Nicas,FMS}}\label{pairs}
Let $(L_1,L_2)$ and $(L'_1,L'_2)$ be two pairs of Lagrangian subspaces of
$\C^n$. Then there exists a unitary map $\pphi \in U(n)$ such that
$\pphi(L_1)=L'_1$ and $\pphi(L_2)=L'_2$ if and only if $\s_{L'_1}\s_{L'_2}$
is conjugate to $\s_{L_1}\s_{L_2}$ in $U(n)$.
\end{thm}

The following series of results will be useful to us in the proof of theorem
\ref{s0-Lag}. The underlying idea is that the elements of the symmetric space
$\lag(n)=U(n)/O(n)$ can be identified with the \emph{symmetric} elements of
$U(n)$ (that is, elements of $U(n)$ satisfying $\tau(u)=u^{-1}$, see
\cite{He,Lo}), all of them being of the form $\pphi^t\pphi$, where $\pphi\in U(n)$
and $\pphi^t$ denotes the transpose of $\pphi$ (so that the symmetric
elements of $U(n)$ are indeed symmetric unitary matrices).

\begin{prop}\label{angles}
Let $W(n):=\{w \in U(n) ~|~ w^t =w\}$ be the set of symmetric unitary matrices.
\begin{enumerate}
\item[(i)] Let $u \in U(n)$. Then $u\in W(n)$ if and only if there exists $k \in
O(n)$ such that $kuk^{-1}$ is diagonal.
\item[(ii)] If $w \in W(n)$, then there exists $\pphi\in W(n)$ such that
$\pphi^2=w$.
\item[(iii)] For any $w\in W(n)$, define $L_w:=\{z \in \C^n ~|~
z-w\overline{z}=0\}$. Then, if $\pphi$ is any element in $W(n)$ such that
$\pphi^2=w$, we have $\pphi(L_0)=L_w$. Consequently, $L_w$ is a Lagrangian
subspace of $\C^n$. Furthermore, $\s_{L_w}\s_{L_0}=w$.
\item[(iv)] The map $w\in W(n) \mapsto L_w \in \lag(n)$ is a diffeomorphism
whose inverse is the well-defined map 
$$\begin{array}{rcl}
\lag(n)=\raisebox{1.25pt}{$U(n)$}\big/\raisebox{-2pt}{$O(n)$} & \longto & W(n) \\
L=u(L_0) & \longmapsto & uu^t
\end{array}$$
\item[(v)] For any $L\in \lag(n)$, we have $\s_{L_0}\s_L=v^t v$, where $v$
is any unitary map such that $v(L)=L_0$.
\item[(vi)] For any $u \in U(n)$, there exist two Lagrangian subspaces $L_1,L_2
\in \lag(n)$ such that $u=\s_{L_1}\s_{L_2}$.
\end{enumerate}
\end{prop}

\begin{proof}
\begin{enumerate}
\item[(i)] Observe that, alternatively, $W(n)=\{w\in U(n) ~|~
w^{-1}=\overline{w}\}$. Now take $w\in W(n)$ and write $w=x+iy$ where $x,y$ are real
matrices. Then $w^t=w$ implies $x^t=x$ and $y^t=y$, and $w\overline{w}=Id$
implies $x^2+y^2=Id$ and $xy-yx=0$. Thus $x$ and $y$ are commuting real
symmetric matrices, so there exists $k \in O(n)$ such that $d_x:=kxk^{-1}$
and $d_y=kyk^{-1}$ are both diagonal. Therefore, $kwk^{-1}=d_x +id_y$ is
diagonal. The converse is obvious. One may observe that since
$d_x^2+d_y^2=k(x^2+y^2)k^{-1}=Id$, one has $d_x+id_y=\exp(iS)$ where $S$ is a
real symmetric matrix.
\item[(ii)] is an immediate consequence of (i).
\item[(iii)] Take $\pphi \in W(n) ~|~ \pphi^2 =w$. Then $z-w\overline{z}=0$ iff
$z-\pphi^2\overline{z}=0$ , that is, $\pphi^{-1}z-\pphi\overline{z}=0$. But
$\pphi^{-1}=\overline{\pphi}$ so that $z\in L_w$ is equivalent to
$\pphi^{-1}z=\overline{\pphi^{-1}z}$ hence to $\pphi^{-1}z\in L_0$, hence to
$z\in\pphi(L_0)$, which shows that $L_w=\pphi(L_0)$ is a Lagrangian subspace
of $\C^n$. Furthermore, $\s_{L_w}\s_{L_0}=\pphi\s_{L_0}\pphi^{-1}\s_{L_0}$.
But since $\s_{L_0}$ is complex conjugation in $\C^n$ and since $\pphi$ is
both symmetric and unitary, we have $\pphi^{-1}\s_{L_0} =
\overline{\pphi^t}\s_{L_0} = (\s_{L_0}\pphi^t\s_{L_0})\s_{L_0}=\s_{L_0}\pphi$,
therefore $\s_{L_w}\s_{L_0}=\pphi\s_{L_0}^2\pphi=\pphi^2=w$.
\item[(iv)] Observe that if $u,v$ are two unitary maps sending $L_0$ to
$L\in\lag(n)$ then $v^{-1}u\in Stab(L_0)=O(n)$ so that $uu^t=vv^t$. Then, if
$L=u(L_0)\in\lag(n)$, one has $L_{uu^t}=\{z-uu^t\overline{z}=0\}$. But
$z-uu^t\overline{z}=0$ iff $u^{-1}z=\overline{u}^{-1}\overline{z}$ , that is, $u^{-1}z \in
L_0$ so $L_{uu^t}=u(L_0)$. Conversely, we know that $L_w=\pphi(L_0)$ where
$\pphi\in W(n) ~|~ \pphi^2=w$ so that indeed $\pphi\pphi^t=\pphi^2=w$.
\item[(v)]For a given $L\in\lag(n)$, take $v\in U(n)$ such that
$v(L)=L_0$. Then $L=v^{-1}(L_0)$ and so we know from (iii) and (iv) that
$L=\{z-(v^{-1})(v^{-1})^t\overline{z}=0\}$ and that
$\s_L\s_{L_0}=v^{-1}(v^{-1})^t$. Hence $\s_{L_0}\s_{L}=(\s_L\s_{L_0})^{-1}=v^tv$.
\item[(vi)] Let $d=\mathrm{diag}\,(\alpha_1,\,\ldots\, ,\alpha_l)\in U(n)$
be a diagonal matrix such that $u=\pphi d^2 \pphi^{-1}$ and set $L=d(L_0)$.
Then we know from (iii) and (iv) that $\s_L\s_{L_0}=d^2$, hence
$u=\pphi\s_{L}\s_{L_0}\pphi^{-1}=\s_{\pphi(L)}\s_{\pphi(L_0)}$.
\end{enumerate}
\end{proof}

Statement (v) may seem a bit useless at this point as it is just a way of
rephrasing (ii), but it will prove useful to us when formulating the
\emph{centered Lagrangian problem} (see subsection \ref{centered}).

\section{Quasi-Hamiltonian spaces}\label{q-Ham}

We recall here the definition of quasi-Hamiltonian spaces and the examples
that shall be useful to us in the following. We follow \cite{AMM} (see also
\cite{GHJW} and \cite{AKSM} for related constructions). Let $U$
\glossary{$U$~: a compact connected Lie group}
be a compact connected Lie group
acting on a manifold $M$ endowed with a $2$-form $\w$. We denote by
$(.\,|\,.)$ an $Ad$-invariant Euclidean product on
$\mathfrak{u}=Lie(U)=T_1\,U$ \glossary{$\mathfrak{u} = Lie(U) = T_1 U$~ the
Lie algebra of $U$} \glossary{$(.|.)$~: an $Ad$-invariant Euclidean scalar
product on $\mathfrak{u}$}
. Let $\chi$ \glossary{$\chi$~: the Cartan $3$-form of $U$}
be (half) the Cartan $3$-form of
$U$, that is, the left-invariant $3$-form defined on $\mathfrak{u}=T_1\,U$ by
$\chi_1(X,Y,Z)=\tfrac{1}{2}\, (X\,|\,[Y,Z])=\tfrac{1}{2}\,([X,Y]\,|\,Z)$,
where the last equality follows from the $Ad$-invariance property. Since
$(.\,|\,.)$ is $Ad$-invariant, $\chi$ is actually bi-invariant and therefore
closed~: $d\chi =0$. Further, denote by $\theta^L$ and $\theta^R$
\glossary{$\theta^L,\theta^R$~: the Maurer-Cartan $1$ forms of $U$}
the left
and right-invariant Maurer-Cartan $1$-forms on $U$~: they take values in
$\mathfrak{u}$ and are the identity on $\mathfrak{u}$, meaning that for any
$u\in U$ and any $\xi \in T_u\,U$, $\theta^L_u(\xi)=u^{-1}.\,\xi$ and
$\theta^R_u(\xi)=\xi.u^{-1}$ (where we denote by a point $.$ the
effect of translations on tangent vectors). Finally, denote by $X^{\sharp}$
\glossary{$X^{\#}$~: the fundamental vector field on a given $U$-space $M$
associated to $X\in\mathfrak{u}$}
the
fundamental vector field on $M$ defined, for any $X\in \mathfrak{u}$, by the
action of $U$~: $X^{\sharp}_x=\tfrac{d}{dt}|_{t=0}(\exp(tX).x)$ for any $x\in
M$. Throughout this paper, we will follow the conventions in \cite{Morita}
to compute exterior products and exterior differentials of differential forms.

\begin{defi}[Quasi-Hamiltonian space]\cite{AMM} In the above notations, $(M,\w)$ is called a
\emph{quasi-Hamiltonian space} if~:
\begin{enumerate}
\item[(i)] The $2$-form $\w$ is $U$-invariant~: $\forall u \in U$, the
associated diffeomorphism of $M$, denoted by $\pphi_u$, satisfies
$\pphi_u^*\w=\w$.
pl\item[(ii)] There exists a map $\mu:M\to U$, called the \emph{momentum map}, such that~:
\begin{enumerate}
\item $\mu$ is equivariant with respect to the $U$-action on $M$ and
conjugation in $U$
\item $d\w=-\mu^*\chi$
\item $\forall x \in M$, $\ker \w_x =\{X^{\sharp}_x \,:\, X\in \mathfrak{u}~|~
Ad\,\mu(x).X=-X\}$
\item $\forall X\in \mathfrak{u}$, the interior product of $X^{\sharp}$ and
$\w$ is $$\iota_{X^{\sharp}}\w = \frac{1}{2}\, \mu^*(\theta^L +
\theta^R\,|\,X)$$ where $(\theta^L+\theta^R\,|\,X)$ is the real-valued
$1$-form defined on $U$ by $(\theta^L +\theta^R\,|\,X)_u(\xi)=
(\theta^L_u(\xi) + \theta^R_u(\xi)\,|\,X)$ for any $u\in U$ and any $\xi \in
T_u\,U$.
\end{enumerate}
\end{enumerate}
\end{defi}

The examples of quasi-Hamiltonian space that will be of most interest to us
are the conjugacy classes of $U$.

\begin{prop}\emph{\cite{AMM}}
Let $\calC\subset U$ be a conjugacy class of a compact connected Lie group
$U$. The tangent space to $\calC$ at $u\in \calC$ is
$T_u\,\calC=\{X.u-u.X\,:\,X\in\mathfrak{u}\}$ . For a given $X\in\mathfrak{u}$,
denote $[X]_u:=X.u-u.X$. Then the $2$-form $\w$ on $\calC$ given at $u\in\calC$
by $$\w_u([X]_u,[Y]_u)= \frac{1}{2}\,
\big((Ad\,u.X\,|\,Y)-(Ad\,u.Y\,|\,X)\big)$$ is well-defined and makes $\calC$ a
quasi-Hamiltonian space for the conjugation action and with momentum map the
inclusion $\mu:\calC\hookrightarrow U$. Such a $2$-form is actually unique.
\end{prop}

Observe that $[X]_u=X^{\sharp}_u$\,, that is~: the fundamental vector fields
generate the tangent space to $\calC$. It is also useful to write this
quantity $[X]_u=(X-Ad\,u.X).u=u.(Ad\,u^{-1}.X - X)$. In order to describe
the symplectic structure on the moduli space $\Mod=\Hom_{\calC}(\pi,U(n))/U(n)$, we
will have to consider the product space $\pconj$ where the $\calC_j$ are
conjugacy classes in $U(n)$, endowed with the diagonal action of $U(n)$. To
make this a quasi-Hamiltonian space with momentum map the product map
$\mu(u_1, \,\ldots\, , u_l)=u_1\ldots u_l$, one has to endow it with a form that is
\emph{not} the product form but has extra terms. The product space thus obtained
is called the \emph{fusion product} and usually denoted
$\calC_1\circledast\cdots\circledast\calC_l$. The general result is the following~:

\begin{thm}[Fusion product of quasi-Hamiltonian spaces]\emph{\cite{AMM}}\label{fusion}
Let $(M_1,\w_1,\mu_1)$ and $(M_2,\w_2,\mu_2)$ be two quasi-Hamiltonian
$U$-spaces. Endow $M_1 \times M_2$ with the diagonal action of $U$. Then the
$2$-form $$\w:=(\w_1\oplus \w_2) + (\mu_1^*\theta^L \wedge
\mu_2^*\theta^R)$$ makes $M_1\times M_2$ a quasi-Hamiltonian space with
momentum map
$$\begin{array}{rccl}
\mu_1\cdot\mu_2: & M_1 \times M_2 & \longto & U \\
& (x_1,x_2) & \longmapsto & \mu_1(x_1)\mu_2(x_2)
\end{array}$$
\end{thm}

\noindent Here, the $2$-form $\w_1\oplus\w_2$ is the product form
$(\w_1\oplus\w_2)_{(x_1,x_2)}((v_1,v_2),(w_1,w_2))= (\w_1)_{x_1}(v_1,w_1) +
(\w_2)_{x_2}(v_2,w_2)$ and $(\mu_1^*\theta^L \wedge \mu_2^*\theta^R)$ is the
$2$-form defined on $M_1\times M_2$ by $$(\mu_1^*\theta^L \wedge
\mu_2^*\theta^R)_{(x_1,x_2)}\big((v_1,v_2),(w_1,w_2)\big) = \frac{1}{2}\,
\Big(\big((\mu_1^*\theta^L)_{x_1}.v_1 \,|\, (\mu_2^*\theta^R)_{x_2}.w_2\big) -
\big((\mu_1^*\theta^L)_{x_1}.w_1 \,|\, (\mu_2^*\theta^R)_{x_2}.v_2\big)\Big)$$

\noindent The above result shows that $\pconj$ is indeed a quasi-Hamiltonian
space for the diagonal action of $U(n)$, with momentum map the product
$\mu(u_1, \,\ldots\, , u_l)=u_1 \ldots u_l$. For a product of three factors,
one can write down the
fusion product form explicitly in the following way~:

\begin{cor}\label{l=3}
The fusion product form on $M_1\times M_2\times M_3$ is the $2$-form
$$\w=(\w_1\oplus\w_2\oplus\w_3) + \big((\mu_1^*\theta^L\wedge\mu_2^*\theta^R)
\oplus (\mu_2^*\theta^L\wedge\mu_3^*\theta^R) \oplus
(\mu_1^*\theta^L\wedge(\mu_2^*Ad).\mu_3^*\theta^R)\big)$$
\end{cor}

\begin{proof}
To obtain the above expression, one applies theorem \ref{fusion}
successively to $M_1\times M_2$ and to $(M_1\times M_2)\times M_3$. One can
then also check that the fusion product is associative, as shown in
\cite{AMM}~: 
\begin{eqnarray*}
\w & = &
\Big(\big((\w_1\oplus\w_2)+(\mu_1^*\theta^L\wedge\mu_2^*\theta^R)\big)\oplus\w_3\Big)
+ \big((\mu_1\cdot\mu_2)^*\theta^L\wedge\mu_3^*\theta^R\big) \\
& = & \Big(\w_1 \oplus
\big((\w_2\oplus\w_3)+(\mu_2^*\theta^L\wedge\mu_3^*\theta^R)\big)\Big)
+ \big(\mu_1^*\theta^L\wedge(\mu_2\cdot\mu_3)^*\theta^R\big)
\end{eqnarray*}
\end{proof}

\section{The symplectic structure on the moduli space of unitary
representations of surface groups}\label{moduli}

The theory of quasi-Hamiltonian spaces provides a very nice description of
the symplectic structure of moduli spaces of unitary representations of
surface groups. We refer to \cite{AMM} for the general description of these
moduli spaces as quasi-Hamiltonian quotients and we will now concentrate on
the space of representations of $\pi=\piS =\, <g_1,g_2,\,\ldots\,
,g_l~|~g_1g_2{\ldots}g_l=1>$. Giving such a representation
with prescribed conjugacy classes $\calC_1, \,\ldots\, , \calC_l$ of
generators amounts to giving $l$ unitary matrices $u_1, \,\ldots\, , u_l$ such
that $u_j \in \calC_j$ and $u_1\ldots u_l=1$. But we know from section
\ref{q-Ham} that this amounts to saying that $(u_1, \,\ldots\, , u_l)\in
\mu^{-1}(\{1\})$ where
$$\begin{array}{rccl}
\mu: & \pconj & \longto & U(n) \\
& (u_1, \,\ldots\, , u_l) & \longmapsto & u_1 \ldots u_l
\end{array}$$
\noindent is the momentum map of the diagonal $U(n)$-action. The moduli
space of unitary representations is then
$\Mod=\Hom_{\calC}(\pi,U(n))/U(n)=\mu^{-1}(\{1\})/U(n)$, which is the symplectic manifold obtained
from $\pconj$ by \emph{quasi-Hamiltonian reduction}, a procedure which we
now recall, stating theorem 5.1 from \cite{AMM} in a particular case
to apply it more directly to our setting.

\begin{thm}[Symplectic reduction of quasi-Hamiltonian
manifolds]\emph{\cite{AMM}}\label{reduction}
Let $(M,\w)$ be a quasi-Hamiltonian $U$-space with momentum
map $\mu:M\to U$. Let
$i:\mu^{-1}(\{1\})\hookrightarrow M$ be the inclusion of the level set
$\mu^{-1}(\{1\})$ in $M$ and let $p:\mu^{-1}(\{1\})\to \mu^{-1}(\{1\})/U$ be
the projection on the orbit space. Assume that $U$ acts freely on $\mu^{-1}(\{1\})$. Then there exists a unique symplectic
form $\w^{red}$ on the reduced space
$M^{red}:=\mu^{-1}(\{1\})/U$ such that $p^*\w^{red}=i^*\w$ on $\mu^{-1}(\{1\})$.
\end{thm}

The proof consists in showing that $i^*\w$ is basic with respect to the
fibration $p$ and then verifying that the corresponding form $\w^{red}$ on
$\mu^{-1}(\{1\})/U$ is indeed a symplectic form. In virtue of the above
theorem, describing the symplectic structure of $\Mod=\mu^{-1}(\{1\})/U$
amounts to giving the $2$-form defining the quasi-Hamiltonian structure on the
product $\pconj$. We now give the description of this $2$-form in the case
where $l=3$.

\begin{prop}\label{2-form}
Let $(u_1,u_2,u_3)\in\calC_1\times\calC_2\times\calC_3$. Take $X_j,Y_j \in
\mathfrak{u}$ and write $[X_j],[Y_j]\in T_{u_j}\,\calC_j$ for the
corresponding tangent vectors (see section \ref{q-Ham}). The $2$-form $\w$
making $\calC_1\times\calC_2\times\calC_3$ a quasi-Hamiltonian space with
momentum map $\mu(u_1,u_2,u_3)=u_1u_2u_3$ is given by~:
\begin{eqnarray*}
\w_u([X],[Y]) & = & \frac{1}{2} \, \big( 
(Ad\,u_1.X_1\,|\,Y_1) - (Ad\,u_1.Y_1\,|\,X_1)
+ (Ad\,u_2.X_2\,|\,Y_2) - (Ad\,u_2.Y_2\,|\,X_2) \\
& & \quad + (Ad\,u_3.X_3\,|\,Y_3) - (Ad\,u_3.Y_3\,|\,X_3)
+ (Ad\,u_1^{-1}.X_1 - X_1 \,|\, Y_2 - Ad\,u_2.Y_2) \\ 
& & \quad - (Ad\,u_1^{-1}.Y_1 - Y_1 \,|\, X_2 - Ad\,u_2.X_2)
+ (Ad\,u_2^{-1}.X_2 - X_2 \,|\, Y_3 - Ad\,u_3.Y_3) \\ 
& & \quad - (Ad\,u_2^{-1}.Y_2 - Y_2 \,|\, X_3 - Ad\,u_3.X_3)
+ (Ad\,u_1^{-1}.X_1 - X_1 \,|\, Ad\,u_2.Y_3 - Ad\,(u_2u_3).Y_3) \\
& & \quad - (Ad\,u_1^{-1}.Y_1 - Y_1 \,|\, Ad\,u_2.X_3 - Ad\,(u_2u_3).X_3)
\big)
\end{eqnarray*}
\end{prop}

The above expression is obtained by applying corollary \ref{l=3}. Observe
that the fusion product $2$-form on $\calC_1\times\calC_2$ consists exactly
of terms of the above expression which do not contain vectors $X_3$ or
$Y_3$. See also remark 5.3 in \cite{Treloar} for expressions of fusion product forms on
products of conjugacy classes. 

\section{Lagrangian submanifolds of a quasi-Hamiltonian
quotient}\label{Lag_QH}

The purpose of this section is to give a way of finding Lagrangian
submanifolds in a symplectic manifold obtained by reduction from a
quasi-Hamiltonian space. It mainly consists in carrying over a standard
procedure for usual symplectic quotients to the quasi-Hamiltonian setting.
To that end, we recall the following result from \cite{OSS} (proposition
2.3), which concerns Hamiltonian spaces. Let $U$ be a compact connected Lie group acting on a symplectic
manifold $(M,\w)$ in a Hamiltonian fashion with momentum map
$\Phi:M\to \mathfrak{u}^*$. Let $\tau$ denote an involutive automorphism of
$U$ and still denote by $\tau$ the involution

$$\begin{array}{rccl}
(T_1\,\tau)^*: & \mathfrak{u}^* & \longto & \mathfrak{u}^* \\
& \lambda & \longmapsto & \lambda \circ T_1\,\tau
\end{array}$$

\noindent that it induces on the dual $\mathfrak{u}^*$ of the Lie algebra
$\mathfrak{u}=T_1\,U$ of $U$. Let $\beta$ be an anti-symplectic involution on
$M$ (that is, such that $\beta^*\w=-\w$ and $\beta^2=Id_M$). In the above
notations, $\beta$ is said to be \emph{compatible with the action} of $U$ if
$\forall u \in U$, $\forall x\in M$, $\beta(u.x)=\tau(u).\beta(x)$ and
$\beta$ is said to be \emph{compatible with the momentum map}
$\Phi:M \to \mathfrak{u}^*$ if $\forall x\in M$, $\Phi\circ\beta(x) =
-\tau\circ\Phi(x)$.

\begin{prop}\emph{\cite{OSS}}\label{fixed_pts}
If $M^{\beta}:=Fix(\beta)$ is non-empty, it is a Lagrangian submanifold of
$M$, stable by the action of the subgroup $U^{\tau}:=Fix(\tau)$ of $U$.
\end{prop}

O'Shea and Sjamaar then proceed to studying the reduced space
$M^{red}=\Phi^{-1}(\{0\})/U$, on which $\beta$ induces an involution
$\hat{\beta}$. To obtain analogous results for
a symplectic manifold $M^{red}=\mu^{-1}(\{1\})/U$ obtained by reduction of a
\emph{quasi}-Hamiltonian space $M$, we wish to define an involution $\beta$ on $M$
such that $\beta$ induces an anti-symplectic involution $\hat{\beta}$ on
$M^{red}$. This is done the following way~:

\begin{thm}\label{lag_locus}
Let $U$ be a compact connected Lie group and let $(M,\w)$ be a
quasi-Hamiltonian $U$-space with momentum map $\mu: M \to U$ . Let
$\tau$ be an involutive automorphism of $U$, denote by $\taum$ the involution
defined on $U$ by $\taum(u)=\tau(u^{-1})$ and let $\beta$ be an involution on
$M$ such that~: 
\begin{enumerate}
\item[(i)] $\forall u \in U$, $\forall x \in M$,
$\beta(u.x)=\tau(u).\beta(x)$ \hfill{\emph{($\beta$ is said to be} compatible
with the action \emph{of $U$)}}
\item[(ii)] $\forall x \in M$, $\mu \circ \beta (x)=\taum \circ \mu (x)$
\hfill{\emph{($\beta$ is said to be} compatible with the momentum map
$\mu:M\to U$\emph{)}}
\item[(iii)] $\beta^* \w = -\w$ \hfill{\emph{($\beta$ reverses the $2$-form
$\w$)}}
\end{enumerate}
then $\beta$ induces an anti-symplectic involution $\hat{\beta}$ on the
reduced space $M^{red}:=\mu^{-1}(\{1\})/U$ . If $\hat{\beta}$ has fixed
points, then $Fix(\hat{\beta})$ is a Lagrangian submanifold of $M^{red}$.
\end{thm}

\begin{rk} See the end of this section for comments on the condition
$Fix(\hat{\beta})\not=\emptyset$.
\end{rk}

\begin{proof}
Compatibility with the momentum map (condition (ii)) shows that $\beta$ maps
$\mu^{-1}(\{1\})$ into $\mu^{-1}(\{1\})$ (since $\taum(1)=1$). Compatibility
with the action (condition (i)) then shows that $\beta(u.x)$ and $\beta(x)$
lie in the same $U$-orbit, so that we have a map

$$\begin{array}{rccl}
\hat{\beta}: & \mu^{-1}(\{1\})/U & \longto & \mu^{-1}(\{1\})/U \\
& U.x & \longmapsto & U.\beta(x)
\end{array}$$

\noindent We know from quasi-Hamiltonian reduction (see theorem \ref{reduction}) that
there exists a unique symplectic form $\w^{red}$ on
$M^{red}=\mu^{-1}(\{1\})/U$ such that $p^*\w^{red}=i^*\w$ where
$i:\mu^{-1}(\{1\})\hookrightarrow M$ and $p:\mu^{-1}(\{1\})\to M^{red}$. To show
that $\hat{\beta}^*\w^{red}=-\w^{red}$, let us first prove that
$i^*(\beta^*\w)$ is basic with respect to the fibration $p$. Then there will
exist a unique $2$-form $\gamma$ on $M^{red}$ such that
$p^*\gamma=i^*(\beta^*\w)$. Since both $\gamma=-\w^{red}$ and
$\gamma=\hat{\beta}^*\w^{red}$ satisfy this condition, they have to be equal.
The last part of the theorem then follows from proposition
\ref{fixed_pts}, as the fixed-point set of an anti-symplectic involution, if
it is non-empty, is always a Lagrangian submanifold. Let us now write this
explicitly.\\ Verifying that $i^*(\beta^*\w)$ is basic is easy since
$\beta^*\w=-\w$ and $i^*\w$ is basic (see \cite{AMM}) but it is actually
true without this assumption so we prove it for $\beta$ satisfying only
conditions (i) and (ii) above. We have to show that $i^*(\beta^*\w)$ is
$U$-invariant and that for every $X\in\mathfrak{u}=Lie(U)$, we have
$\iota_{X^{\sharp}}(i^*(\beta^*\w))=0$, where $X^{\sharp}$ is as usual the fundamental vector
field $X^{\sharp}_x=\frac{d}{dt}|_{t=0}(\exp(tX).x)$ (for any $x\in M$) associated
to $X\in\mathfrak{u}$ by the action of $U$ on $M$. Let $u\in U$ and denote
by $\pphi_u$ the corresponding diffeomorphism of $M$. The map $\mu$ being
equivariant $\pphi_u$ sends $\mu^{-1}(\{1\})$ into itself, hence
$i\circ\pphi_u=\pphi_u\circ i$ on $\mu^{-1}(\{1\})$. Furthermore,
compatibility with the action yields $\beta \circ \pphi_u =
\pphi_{\tau(u)}\circ \beta$. We then have, on $\mu^{-1}(\{1\})$, 
\begin{eqnarray*}
\pphi_u^*(i^*(\beta^*\w)) & = & (\beta \circ i \circ \pphi_u)^* \w \\
& = & (\pphi_{\tau(u)} \circ \beta \circ i)^* \w \\
& = & i^*(\beta^*(\underbrace{\pphi_{\tau(u)}^*\w}_{=\w}))
\end{eqnarray*}
\noindent where the very last equality follows from the $U$-invariance of
$\w$. Further, let $X\in\mathfrak{u}$. Since $\beta$ is compatible with the
action, one has
$\beta(\exp(tX).x)=\tau(\exp(tX)).\beta(x)=\exp(t\tau(X)).\beta(x)$ (where we
still denote by $\tau$ the involution $T_1\,\tau$ on $\mathfrak{u}=T_1\,U$),
hence $T_x\beta.X^{\sharp}_x=(\tau(X))^{\sharp}_{\beta(x)}$, hence
$\iota_{X^{\sharp}}(\beta^*\w) = \beta^*(\iota_{(\tau(X))^{\sharp}}\w)$. Since
$\iota_{X^{\sharp}}(i^*(\beta^*\w))= i^*(\iota_{X^{\sharp}}(\beta^*\w))$, we
can compute, using the fact that $\beta$ is compatible with $\mu$,

\begin{eqnarray*}
\iota_{X^{\sharp}}(\beta^*\w) & = & \beta^*(\iota_{(\tau(X))^{\sharp}}\w) \\
& = & \beta^*\Big(\frac{1}{2}\,\mu^*\big(\theta^L + \theta^R \, | \, \tau(X)\big)\Big) \\
& = & \frac{1}{2}\, (\mu\circ\beta)^*\big(\theta^L + \theta^R \,|\,
\tau(X)\big) \\
& = & \frac{1}{2}\, (\taum\circ\mu)^*\big(\theta^L + \theta^R \,|\,
\tau(X)\big) \\
& = & \frac{1}{2}\, \mu^* \Big((\taum)^*\big(\theta^L + \theta^R \,|\,
\tau(X)\big)\Big)
\end{eqnarray*}

\noindent hence $i^*(\iota_{X^{\sharp}}(\beta^*\w)) = \frac{1}{2}\,
i^*\circ\mu^*\, (\ldots) = \frac{1}{2}\, (\mu\circ i)^*\,(\ldots)$. But
$\mu\circ i : \mu^{-1}(\{1\})\to U$ is a constant map, therefore $T(\mu\circ
i)$ and consequently $(\mu\circ i)^*$ are zero, which completes the proof
that $i^*(\beta^*\w)$ is basic. Finally, let us show that
$p^*(\hat{\beta}\w^{red}) = i^*(\beta^*\w) = p^*(-\w^{red})$ (this is where
we really use $\beta^*\w=-\w$). We have, on $\mu^{-1}(\{1\})$,
$p^*(\hat{\beta}^*\w^{red})=(\hat{\beta}\circ p)^*\w^{red}=(p\circ
\beta)^*w^{red}= \beta^* (p^*\w^{red})=\beta^*(i^*\w)=(i\circ\beta)^* \w =
(\beta\circ i)^*\w = i^*(\beta^*\w)=i^*(-\w)=-i^*\w=-p^*\w^{red} =
p^*(-\w^{red})$. This completes the proof, as indicated above.
\end{proof}

In the next section, we will give an example of a map $\beta$ satisfying the
hypotheses of theorem \ref{lag_locus}. In the case we will then be dealing
with, it will be important to us that the map $\tau$ under consideration be
actually an isometry for the Euclidean product $(.\,|\,.)$ on $\mathfrak{u}$
(recall that this scalar product is part of the initial data to define
quasi-Hamiltonian $U$-spaces). So far, we did not need that hypothesis.
Before ending this section, we would like to say that, in fact,
if $\beta$ satisfies the conditions of theorem \ref{lag_locus} and has
fixed points then $\hat{\beta}$ necessarily has fixed points. Indeed,
observe first that $Fix(\hat{\beta})\not=\emptyset$ if and only if
$Fix(\beta)\cap\mu^{-1}(\{1\})\not=\emptyset$. We then have the following
result, which is a convexity result concerning momentum maps in the
quasi-Hamiltonian framework and which is adapted from (a special case of)
the convexity theorem of O'Shea and Sjamaar (see \cite{OSS})~:

\begin{thm}\emph{\cite{mythesis}}\label{convexity}
Let $\beta$ be an involution defined on a quasi-Hamiltonian $(U,\tau)$-space
$(M,\w,\mu:M\to U)$ such that $\beta$ is compatible with the action and the
momentum map and such that $\beta^*\w=-\w$. Assume that
$Fix(\beta)\not=\emptyset$ and that there exists a
maximal torus $T$ of $U$ which is fixed pointwise by $\taum$ and let $\mathcal{W}\subset
\mathfrak{t}=Lie(T)$ be a Weyl alcove. Then $\mu(M^{\beta})\cap \exp\mathcal{W} =
\mu(M)\cap \exp\mathcal{W}$.
\end{thm}
 The proof of this theorem is too long to be presented here, all the more so
as it calls for techniques which are very different from the ones we have
made use of so far. A proof of this result is available in \cite{mythesis}
and will appear in a forthcoming paper (\cite{S}). The fact that
$Fix(\hat{\beta})\not=\emptyset$ is then a corollary of this theorem~:
\begin{cor}\emph{\cite{mythesis}}
If $Fix(\beta)\not=\emptyset$ and $\mu^{-1}(\{1\})\not=\emptyset$ then
$Fix(\beta)\cap\mu^{-1}(\{1\})\not=\emptyset$, in which case the involution
$\hat{\beta}$ induced by $\beta$ on $\mu^{-1}(\{1\})/U$ satisfies
$Fix(\hat{\beta})\not=\emptyset$.
\end{cor}
\begin{proof}[Proof of the corollary]
Since $Fix(\beta)\not=\emptyset$, the above claim applies. Since
$\mu^{-1}(\{1\})\not=\emptyset$ and $1\in \exp\mathcal{W}$, we
then have $1\in \mu(M)\cap \exp\mathcal{W} = \mu(M^{\beta})
\cap \exp\mathcal{W}$, which means that $\mu^{-1}(\{1\})\cap Fix(\beta)
\not=\emptyset$, which in turn is equivalent to $Fix(\hat{\beta})\not=
\emptyset$.
\end{proof}
\noindent We will not use theorem \ref{convexity}, nor its corollary, in the following.

\section{Lagrangian representations as fixed point set of an involution}

In this section, we will state our main result, which is the
characterization of a $\s_0$-Lagrangian representations of $\pi=\piS$ as the
elements of the fixed point set of an involution $\beta$ defined on the
product of $l$ conjugacy classes of the unitary group (satisfying the
condition $\exists~ (u_1, \,\ldots\, , u_l)\in \calC_1\times\cdots\times\calC_l
~|~ u_1 \ldots u_l=1$). Using theorem \ref{lag_locus} proved in the
preceding section, we will deduce that the set of (equivalence classes of)
Lagrangian representations is a Lagrangian submanifold of the moduli space
$\Mod=\Hom_{\calC}(\pi,U(n))/U(n)$. The first five subsections explain how these
results (in particular the involution $\beta$) were obtained but may be
skipped if one wants to go straight to the actual theorems (whose proofs may
also be read without knowledge of the previous subsections).

\subsection{The infinitesimal picture and the momentum map approach}

Let us recall our problem~: given $l$ unitary matrices $u_1, \,\ldots\, , u_l
\in U(n)$ satisfying $u_j\in\calC_j$ \glossary{$\calC_1, \,\ldots\, ,
\calC_l$~: $l$ conjugacy classes in $U$}
and $u_1\ldots u_j=1$, do there exist
$l$ Lagrangian subspaces $L_1, \,\ldots\, , L_l$ of $\C^n$ such that
$\s_j\s_{j+1}=u_j$ (where $\s_j$ is the Lagrangian involution associated
with $L_j$ and $\s_{l+1}=\s_1$)\,? As was recalled in section \ref{laginv},
the condition $\s_j\s_{j+1} \in \calC_j$, which lies on the spectrum of
the unitary map $\s_j\s_{j+1}$, can be interpreted geometrically as the
measure of an angle between Lagrangian subspaces. The Lagrangian problem
above can therefore be thought of as a configuration problem in the
Lagrangian Grassmannian $\lag(n)$ of $\C^n$~: given eigenvalues
$\exp(i\lambda_j)$, $\lambda_j \in \R^n$, do there exist $l$ Lagrangian
subspaces $L_1, \,\ldots\, , L_l$ such that
$\mathrm{measure}(L_j,L_{j+1})=\exp(i\lambda_j)$\,? Under this geometrical form,
the Lagrangian problem is slightly more general than our original
representation theory problem. It is very much linked to the unitary problem
studied for instance in \cite{JW,gal,Bi1,AW,KM,Bi2,Be}, which is the following~:
given $\lambda_j\in\R^n$, do there exist $l$ unitary matrices $u_1,
\,\ldots\, , u_l$ satisfying $\mathrm{Spec}\,u_j=\exp(i\lambda_j)$ and
$u_1\ldots u_l =1$~? In fact, a solution $(L_1, \,\ldots\, , L_l)$ to the
Lagrangian problem (second version) provides a solution $u_j=\s_j\s_{j+1}$
to the unitary problem. As was shown in \cite{FMS}, it is possible to use
this approach to give an interpretation of the inequalities found by Biswas
in \cite{Bi1} (which are necessary and sufficient conditions on the
$\lambda_j$ for the unitary problem to have a solution in the case $n=2$ and
$l=3$) in terms of the inequalities satisfied by the angles of a spherical
triangle.\\
The fact that the unitary problem admits a symplectic description was our
first motivation to study the Lagrangian problem from a symplectic point of
view. The second motivation is derived from the above-given geometrical
formulation of the problem. To better understand this, let us try and
formulate an infinitesimal version of the Lagrangian problem. Take three
Lagrangian subspaces $L_1,L_2,L_3$ close enough so that we can think of these points
in $\lag(n)$ as tangent vectors to $\lag(n)$ at some point $L_0$
representing the center of mass of $L_1,L_2,L_3$. Tangent vectors to the
Lagrangian Grassmannian are identified with real symmetric matrices
$S_1,S_2,S_3$ and the center of mass condition then turns into
$S_1+S_2+S_3=0$. It seems reasonable in this context to translate the angle
condition $\mathrm{mes}(L_j,L_{j+1})=\exp(i\lambda_j)$ (that is,
$\mathrm{Spec}\,\s_j\s_{j+1}=\exp(i\lambda_j))$ into the spectral condition
$\mathrm{Spec}\,S_j=\lambda_j\in\R^n$. We then recognize a real version
(replacing complex Hermitian matrices with real symmetric ones) of a famous
problem in mathematics (see \cite{Fu} for a review of this problem and those
related to it)~: given $\lambda_j \in \R^n$, do there exist Hermitian
matrices $H_1,H_2,H_3$ such that $\mathrm{Spec}\,H_j=\lambda_j$ and
$H_1+H_2+H_3=0$~? In fact, these last two problems are equivalent (meaning
that, for given $(\lambda_j)_j$, one of them has a solution if and only if
the other one does) and this can be shown in a purely symplectic framework (see
\cite{AMW}) using momentum maps to translate the condition $H_1+H_2+H_3=0$
into $(H_1,H_2,H_3)\in\mu^{-1}(\{0\})$. Therefrom, it seems promising to try
to think of the Lagrangian problem as a real version, in a sense that will
be made precise in subsection \ref{real_problems}, of the unitary problem
(since a solution to the Lagrangian problem provides an obvious solution to
the unitary problem). 

\subsection{The centered Lagrangian problem}\label{centered}

As a consequence of the above infinitesimal picture, we replace our
Lagrangian problem with a centered problem, meaning that instead of
measuring the angles $(L_j,L_{j+1})$, we measure the angles $(L_0,L_j)$
where $L_0$ is the horizontal Lagrangian $L_0=\R^n\subset\C^n$ (playing the
role of an origin in $\lag(n)$). Recall from section \ref{laginv} (theorem
\ref{pairs} and proposition \ref{angles}) that this
angle is measured by the spectrum of $\s_0\s_j=u_j^tu_j$, where $u_j$ is any
unitary map sending $L_j$ to $L_0$. We then ask the following question~:
given $l$ conjugacy classes $\calC_1, \,\ldots\, , \calC_l \subset U(n)$, does
there exist $l$ unitary matrices $u_1, \,\ldots\, , u_l$ such that $u_j^tu_j
\in \calC_j$ and $u_1\ldots u_l=1$~? The main observation here is then to
see that the condition $\mathrm{Spec}\,u^tu=\exp(i\lambda)$, for some
$\lambda\in\R^n$ (that is, $u^tu$ lies in some fixed conjugacy class of
$U(n)$) means that $u$ belongs to a fixed orbit of the action of $O(n)\times
O(n)$ on $U(n)$ given by $(k_1,k_2).u=k_1uk_2^{-1}$, as is shown by the
following elementary result~:
\begin{lem}\label{dbl_coset}
For any $u,v\in U(n)$, $\mathrm{Spec}\,u^tu=\mathrm{Spec}\,v^tv$ if and only
if there exist $(k_1,k_2)\in O(n)\times O(n)$ such that $v=k_1 u k_2^{-1}$.
\end{lem}
Since we think of the above problem as a real version of some complex
problem, we now wish to find this complex version, which is done by
abstracting a bit our situation to put it in the appropriate framework.

\subsection{Complexification of the centered Lagrangian
problem}\label{complexification}

Let us formulate the centered Lagrangian problem in greater generality. For
everything regarding the theory of Lie groups and symmetric spaces,
especially regarding real forms and duality, we refer to \cite{He}. We start
with a real Lie group $H$. Let $G=H^{\C}$ be its complexification and let
$\tau$ be the Cartan involution on $G$ associated to $H$, that is to say, the
involutive automorphism of $G$ such that $Fix(\tau)=H$. Let $U$ be a compact
connected real form of $G$ such that the associated Cartan involution $\theta$
satisfies $\theta\tau = \tau\theta$. Such a compact group always exists and
is stable under $\tau$. The group $H$ is then stable under $\theta$ and $U$ and
$H$
\glossary{$G=U^{\C}$~: the complexification of $U$}
\glossary{$H$~: the non-compact dual of $U$ in $G$}
\glossary{$\theta$~: the involution on $G$ such that $U=Fix(\theta)$}
\glossary{$\tau$~: the involution on $G$ such that $H=Fix(\tau)$}
\glossary{$K=U\cap H \subset G$}
are said to be dual to each other (when $H$ is non-compact, they indeed
define dual symmetric spaces $U/(U\cap H)$ and $H/(U\cap H)$). Moreover,
because of the fact that $\tau$ is the Cartan involution associated to the non-compact
dual $H$ of $U$, the compact connected group $U$ contains a maximal torus
$T$ such that $\tau(t) = t^{-1}$ for all $t\in T$ ($(U,\tau)$ is said to be
\emph{of maximal rank}, see \cite{Lo} pp. 72-74 and 79-81). Let $K:=U\cap
H$. Then $K=Fix(\tau|_U)\subset U$ and $K=Fix(\theta|_H)\subset H$. We
consider the action of $K\times K$ on $U$ given by $(k_1,k_2).u=k_1 u
k_2^{-1}$. Notice that if $H$ is compact to start with, then $K=U=H$ and the
above action defines congruence in $U$. As for us though, we are interested
in the case where $H$ is non-compact. For $H=Gl(n,\R)$, we have $U=U(n)$ and
$K=O(n)$, and we are then led to ask the following question, which
is a generalized version of our centered Lagrangian problem\,: given $l$
orbits $\calD_1, \,\ldots\, , \calD_l$ \glossary{$\calD_1, \,\ldots\, ,
\calD_l$~: $K\times K$-orbits in $U$} of the action of $K \times K$ on $U$,
do there exist $u_1, \,\ldots\, , u_l \in U$ such that $u_j \in \calD_j$ and
$u_1\ldots u_l=1$~? Observe that, as a generalization of lemma
\ref{dbl_coset}, these orbits are in one-to-one correspondence with the
conjugacy classes in $U$ of elements of the form $\taum(u)u$, where $u$ is
any element in a given orbit $\calD$ and $\taum(u)=\tau(u^{-1})$. Indeed,
this is a corollary of theorem 8.6 in chapter VII (p. 323) of \cite{He},
which we now state under a form most convenient for our purposes~:

\begin{thm}[Cartan decomposition of $U$]\label{Cartan_decomp}\emph{\cite{He}}
Let $U$ be a compact connected Lie group and let $\tau$ be an involutive
automorphism of $U$. Let $K=Fix(\tau)\subset U$. Still denote by $\tau$ the
involutive automorphism $T_1\,\tau:\mathfrak{u}=T_1\,U\to\mathfrak{u}$ . Then
there exist a subset $\mathfrak{q}_0 \subset \mathfrak{u}$ such that~: 
\begin{enumerate}
\item[(i)] $\forall X \in \mathfrak{q}_0$, $\tau(X)=-X$
\item[(ii)] each $u\in U$ can be written $u=k_1\exp(X)k_2^{-1}$ for some
$k_1,k_2 \in K$ and for a unique $X \in \mathfrak{q}_0$.
\end{enumerate}
Further, if $X,Y\in\mathfrak{u}$ satisfy $\tau(X)=-X$ and $\tau(Y)=-Y$, and
if there exist $u\in U$ such that $Ad\,u.X=Y$, then there exists $k\in
K\subset U$ such that $Ad_U\,k.X=Y$.
\end{thm}

\begin{cor}\label{KxK-orbits}
Let $u,v\in U$. Then there exist $(k_1,k_2) \in K\times K$ such that $v=k_1
u k_2^{-1}$ if and only if $\taum(v)v$ and $\taum(u)u$ lie in a same
conjugacy class in $U$.
\end{cor}

\begin{proof}
The first implication is obvious. Conversely, write $u=k_1\exp(X) k_2^{-1}$
as in the above theorem. Then $\taum(u)=k_2\exp(X)k_1^{-1}$ (since
$\tau(k_j)=k_j$ in $U$ and $\taum(X)=-\tau(X)=X$ in $\mathfrak{u}$) and
therefore $\taum(u)u=k_2\exp(2X)k_2^{-1}$. Likewise, we can write
$v=k'_1\exp(Y)(k'_2)^{-1}$ and therefore
$\taum(v)v=k'_2\exp(2Y)(k'_2)^{-1}$. Since $\taum(v)v$ is conjugate to
$\taum(u)u$ in $U$, we see that $2Y$ is $Ad_U\,U$-conjugate to $2X+H$ where
$H\in\mathfrak{u}$ satisfies $\exp(H)=1$. We then necessarily have
$\tau(H)=-H$. By using theorem 8.5 in chapter VII of \cite{He}, we can then
write $H=2Z$ with $Z\in\mathfrak{u}$ satisfying $\tau(Z)=-Z$ and $\exp(Z)\in
K$. Then $Y$ is $Ad_U\,U$-conjugate to $(X+Z)$. But $\tau(Y)=-Y$ and
$\tau(X+Z)=-(X+Z)$, therefore, by the above theorem, $Y$ and $(X+Z)$ are
$Ad_U\,K$-conjugate. Then we have $Y=k.(X+Z).k^{-1}$ in $\mathfrak{u}$ for
some $k\in K$, so that $v=k'_1\exp(Y)(k'_2)^{-1} =
k'_1k\exp(X)\exp(Z)k^{-1}(k'_2)^{-1} = \underbrace{k'_1kk_1^{-1}}_{\in K}
(\underbrace{k_1\exp(X)k_2^{-1}}_{=u})
\underbrace{k_2\exp(Z)k^{-1}(k'_2)^{-1}}_{\in K}$
\end{proof}

Now, to find the complex version of our problem, we apply the same
construction to the complex Lie group $G=H^{\C}$ viewed as a real Lie group.
Then $G^{\C}=G\times G$ is the complexification of $G$ and
$\widetilde{U}=U\times U \subset G \times G = G^{\C}$
\glossary{$\widetilde{U}=U\times U$~: the compact real form of
$G^{\C}=G\times G$}
is a compact real form of
$G^{\C}$. Its non-compact dual (which needs to be a subgroup of $G^{\C}=G
\times G$) is then $\widetilde{H}=\{(g,\theta(g))\,:\,g\in G\}\simeq G$
\glossary{$\widetilde{H}$~: the non-compact dual of $\widetilde{U}$ in $G^{\C}$}
where
$\theta$ is the Cartan involution associated to $U$. The Cartan involution
associated to $\widetilde{U}$ is
$\widetilde{\theta}:(g_1,g_2)\in\widetilde{G} \mapsto
(\theta(g_1),\theta(g_2))$ \glossary{$\widetilde{\theta}$~: the involution
on $G^{\C}$ such that $Fix{\widetilde{\theta}}=\widetilde{U}$}
and the Cartan involution associated to
$\widetilde{H}$ is $\widetilde{\tau}:(g_1,g_2)\mapsto
(\theta(g_2),\theta(g_1))$.\glossary{$\widetilde{\tau}$~: the involution
on $G^{\C}$ such that $Fix{\widetilde{\tau}}=\widetilde{H}$}
Indeed, $Fix(\widetilde{\theta})=\widetilde{U}$,
$Fix(\widetilde{\tau})=\widetilde{H}$ and
$\widetilde{\theta}\widetilde{\tau}=\widetilde{\tau}\widetilde{\theta}$. Then
we define~:
\begin{eqnarray*}
\widetilde{K} & := & \widetilde{U}\cap\widetilde{H} \\
& = & \Big\{\big(g,\theta(g)\big)~|~\widetilde{\theta}\big(g,\theta(g)\big)
= \big(g,\theta(g)\big)\Big\} \\
& = & \Big\{\big(g,\theta(g)\big)~|~\theta(g)=g\Big\} \\
& = & \big\{(u,u)\,:\, u\in U\big\}
\end{eqnarray*}
\glossary{$\widetilde{K} := \widetilde{U}\cap\widetilde{H}$}
(we will also use the notation $U_{\Delta}:=\{(u,u)\,:\, u\in U\}$ instead
of $\widetilde{K}$) and we consider the action of $\widetilde{K}\times\widetilde{K} =
U_{\Delta}\times U_{\Delta}$ on $\widetilde{U}=U\times U$ defined by~:
$$\big((u_1,u_1),(u_2,u_2)\big).(u,v)=(u_1 u u_2^{-1}, u_1 v u_2^{-1})$$
Our problem then states~: given $l$ orbits $\calDt_1, \,\ldots\, , \calDt_l$
\glossary{$\calDt_1, \,\ldots\, , \calDt_l$~:
$\widetilde{K}\times\widetilde{K}$-orbits in $\widetilde{U}$}
of the above action, do there exist $l$ pairs $(u_1,v_1), \,\ldots\, ,$ 
$(u_l,v_l) \in \widetilde{U}=U\times U$ such that $(u_j,v_j)\in\calDt_j$ and
$(u_1,v_1) . \,\ldots\, . (u_l,v_l) = 1$, that is, $u_1 \ldots u_l=1$ and $v_1
\ldots v_l =1$~? \\
Before passing on to the next subsection, we wish to point out that if we
consider the action of $K \times K$ not on $U$ but rather on its dual $H$, then the
orbits of this action are characterized by the singular values
($\mathrm{Sing}\,h=\mathrm{Spec}\,(\theta^{-}(h)h)$ where $h\in H$ and
$\theta^{-}(h)=\theta(h^{-1})$) of any of their elements. As a consequence,
our (centered) Lagrangian problem appears as a compact version of the (real)
Thompson problem, replacing $\theta$ with $\tau$ in the latter to formulate
the former (see \cite{AMW} and \cite{Lu-Ev} for a proof of the
Thompson conjecture in the real case).

\subsection{Equivalence between the complexification of the centered
Lagrangian problem and the unitary problem}

From now on, the initial data is a compact connected Lie group $U$. For such
a group, we can formulate~: (i) the centered Lagrangian problem (concerning
$K \times K$-orbits in $U$, where $K=U\cap H$ with $H$ the non-compact dual
of $U$) (ii) a complex version of this (concerning $U_{\Delta} \times
U_{\Delta}$-orbits in $U\times U$) (iii) the unitary problem (concerning
conjugacy classes in $U$). To show the equivalence of these last two
problems, the main observation to make is the following one~:

\begin{lem}\label{eta-map}
The map
$$\begin{array}{rccl}
\eta : & U \times U & \longto & U \\
& (u,v) & \longmapsto & u^{-1}v
\end{array}$$
sends a $U_{\Delta}\times U_{\Delta}$-orbit $\widetilde{\calD}$ in $U\times U$ onto a
conjugacy class $\calC$ in $U$.
\end{lem}

\begin{proof}
If $(u,v)=(u_1 u_0 u_2^{-1}, u_1 v_0 u_2^{-1})$ then
$u^{-1}v=u_2(u_0^{-1}v_0)u_2^{-1}$ so $\eta(\calDt)\subset \calC$ where
$\calC$ is a conjugacy class in $U$. Further, let $(u,v)\in\calDt$ and take
any $w\in\calC$. Then $\exists\, u_2 ~|~ w=u_2 u^{-1}v u_2^{-1}$ so that
$(1,w)=(1,u_2u^{-1}vu_2^{-1})=\underbrace{((u_2u^{-1},u_2u^{-1}),(u_2,u_2))}_{\in
U_{\Delta}\times U_{\Delta}}.\underbrace{(u,v)}_{\in \calDt}$, hence
$(1,w)\in\calDt$, therefore $w=\eta(1,w)\in\eta(\calDt)$.
\end{proof}

It is nice to observe that this map $\eta$ may be used to show that a
compact connected Lie group $U$ is a symmetric space $U=(U\times
U)/U_{\Delta}$ (see \cite{He}). Coming back to the matter at hand, we have
the following result, that says that the complexification of the centered
Lagrangian problem has a solution if and only if the unitary problem has a
solution (that is, these two problems are equivalent).

\begin{prop}\label{equiv}
Let $\calDt_1, \,\ldots\, , \calDt_l$ be $l$ orbits of $U_{\Delta}\times
U_{\Delta}$ in $U\times U$ and let $\calC_1, \,\ldots\, ,\calC_l \subset U$ be
the corresponding conjugacy classes~: $\calC_j=\eta(\calDt_j)$. Then there
exists $((u_1,v_1), \,\ldots\, , (u_l,v_l)) \in \pdblt$ such that $u_1\ldots
u_l=1$ and $v_1 \ldots v_l=1$ if and only if there exist $(w_1,\,\ldots\,
,w_l) \in \pconj$ such that $w_1\ldots w_l=1$.
\end{prop}

\begin{proof}
Setting $(u_j,v_j):=(1,w_j)$ for every
$j$, we see that the second condition implies the first one. Conversely,
assume that $((u_1,v_1),\,\ldots\, ,(u_l,v_l)) \in \calDt_1 \times \cdots
\times \calDt_l$ satisfy $u_1 \ldots u_l=1$ and $v_1\ldots v_l=1$. Then
$(u_1 \ldots u_l)^{-1}v_1\ldots v_l=1$, hence
$u_l^{-1}\ldots u_2^{-1}(u_1^{-1}v_1)v_2 \ldots v_l=1$, with
$u_1^{-1}v_1\in\calC_1$. Hence~:
$$\underbrace{u_l^{-1}\ldots u_2^{-1}(u_1^{-1}v_1)u_2\ldots u_l}_{\in\calC_1}
\,\underbrace{u_l^{-1}\ldots u_3^{-1}(u_2^{-1}v_2)u_3\ldots
u_l}_{\in\calC_2} \, \ldots \,
\,\underbrace{(u_l^{-1}v_l)}_{\in\calC_l}=1$$
Setting $w_1=u_l^{-1}\ldots u_2^{-1}(u_1^{-1}v_1)u_2\ldots u_l$,
$w_2=u_l^{-1}\ldots u_3^{-1}(u_2^{-1}v_2)u_3\ldots
u_l$, \ldots, and $w_l=u_l^{-1}v_l$ then gives a
solution $(w_1,\,\ldots\, ,w_l)$ to the unitary problem.
\end{proof}

In analogy with a result on double cosets of $U(n)$ in $Gl(n,\C)$ (which are
characterized by the singular values
$\mathrm{Sing}\,g=\mathrm{Spec}\,(\theta^{-}(g)g)$ of any of their elements) and
dressing orbits of $U(n)$ in $(U(n))^*=\{b\in Gl(n,\C)~|~ b~
\mathrm{is~upper~triangular~and~} diag(b) \in (\R^{*+})^n\}$ appearing
in \cite{AMW}, the above proposition can be formulated more precisely in the
following way. Consider the action of $U^l$ on $\pdblt$ given by
$$(\pphi_1,\,\ldots\, ,\pphi_l).\big((u_1,v_1),\,\ldots\,
,(u_l,v_l)\big)=\big(\!\!\!\underbrace{\pphi_1.(u_1,v_1).\pphi_2^{-1}}_{=(\pphi_1 u_1
\pphi_2^{-1}, \pphi_1 v_1 \pphi_2^{-1})},
\pphi_2.(u_2,v_2).\pphi_3^{-1}, \,\ldots\, ,
\,\pphi_l.(u_l,v_l).\pphi_l^{-1}\big)$$ 
and the diagonal action of $U(n)$ on $\pconj$~: $\pphi.(w_1, \,\ldots\, ,
w_l)=(\pphi w_1 \pphi^{-1}, \,\ldots\, , \pphi w_l \pphi^{-1})$. These
actions respectively preserve the relations $u_1 \ldots u_l = v_1 \ldots
v_l=1$ and $\w_1 \ldots \w_l=1$. We may then define the orbit spaces 
$$\mathcal{M}_{\calDt}:=\raisebox{2.5pt}{$\Big\{\big((u_j,v_j)\big)_j\in\pdblt ~|~ u_1\ldots u_l =
v_1\ldots v_l =1\Big\}$}\bigg/\raisebox{-4pt}{$U^l$}$$
and $$\Mod=\raisebox{2pt}{$\big\{(w_j)_j \in \pconj ~|~ w_1\ldots
w_l=1\big\}$}\Big/\raisebox{-3.2pt}{$U$}$$
And we then have

\begin{prop}\label{Thompson1}
The map
$$\begin{array}{rccl}
\eta^{(l)} : & \pdblt & \longto & \pconj \\
& ((u_1,v_1), \,\ldots\, , (u_l,v_l)) & \longmapsto & (u_l^{-1}\ldots
u_2^{-1}(u_1^{-1}v_1)u_2\ldots u_l, \,\ldots\, ,
\, u_l^{-1}(u_{l-1}^{-1}v_{l-1})u_l,\, u_l^{-1}v_l)
\end{array}$$
induces a homeomorphism $\mathcal{M}_{\calDt}\simeq\Mod$.
\end{prop}
We will not use this result in the following so we do not give the proof,
which is but a consequence of the above. We point out the fact that this
result reinforces the analogy between our problem and the Thompson problem.
We now wish to explain in what precise sense the Lagrangian problem is a
real version of these two equivalent problems.

\subsection{Solutions to real problems as fixed point sets of
involutions}\label{real_problems}

The important idea of thinking of possible solutions to a real problem
as the fixed point set of an involution defined on the set of possible
solutions to a corresponding complex problem is well-established in
symplectic geometry and is due to Michael Atiyah and Alan Weinstein (see
\cite{Ati,Du} and \cite{Lu-Ra}).
In fact, the idea is that the set of possible solutions to a complex problem
carries a symplectic structure and that the corresponding real problem is
formulated for elements of the fixed point set of an anti-symplectic
involution defined on this symplectic manifold. Examples of results obtained
using this idea include the (linear and non-linear) real Kostant convexity theorems
(see \cite{Du,Lu-Ra}) and the real Thompson conjecture (see
\cite{AMW,Lu-Ev}). Although we will have to replace symplectic manifolds with
quasi-Hamiltonian spaces for technical considerations, the above idea plays
a key role in our approach. Keeping this in mind, we will eventually define
an involution $\beta^{(l)}$ on the quasi-Hamiltonian space $\pconj$. But, to
explain how this involution is obtained, we will first work on the product
$\pdblt$ of $l$ $U_{\Delta}\times U_{\Delta}$-orbits in $U\times U$.\\
The key here is to try and see the $K\times K$-orbit of $w\in U$ as a subset
of some $U_{\Delta}\times U_{\Delta}$-orbit $\calDt \subset U \times U$. This is done
by observing that $w\in\calD$ is equivalent to $\taum(w)w\in\calC$ which in
turn is equivalent to $(\tau(w),w)\in\calDt$. Indeed, the first equivalence
is corollary \ref{KxK-orbits}, where $\calC$ is defined as the conjugacy class
of $\taum(w)w$ for any $w\in\calD$. Then we know from lemma \ref{eta-map}
that $\calC=\eta(\calDt)$ where $\calDt$ is the $U_{\Delta}\times
U_{\Delta}$-orbit of $(1,\taum(w)w)\sim(\tau(w),w)$, which gives the second
equivalence. In order to obtain elements of the form $(\tau(w),w)$ as fixed
points of an involution, we set

$$\begin{array}{rccl}
\alpha: & U\times U & \longto & U \times U \\
& (u,v) & \longmapsto & \big(\tau(v),\tau(u)\big)
\end{array}$$

\noindent Then $\alpha^2=Id$ and $Fix(\alpha)=\{(\tau(v),v)~|~v\in U\}\simeq U$. In
particular, $Fix(\alpha)$ is always non-empty. Moreover, we have~:

\begin{lem}\label{D-inv}
$\alpha(\calDt)=\calDt$, so that $\alpha$ defines an involution on $\calDt$,
whose fixed point set is isomorphic to $\calD$ and therefore non-empty.
\end{lem}

\begin{proof}
If $(u,v)\in\calDt$, we have $\eta(\alpha(u,v))=\taum(v)\tau(u)=
\tau(v^{-1}u) =\taum(u^{-1}v)$. But if $w\in U$, then $\taum(w)$ is
conjugate to $w$. Indeed, since $\tau$ comes from the Cartan involution defining the non-compact dual of $U$, there always exists a
maximal torus of $U$ which is fixed pointwise by $\taum$ (see \cite{Lo}, pp. 72-74 and 79-80),
and $w$ is conjugate to
an element in such a torus~: $w=\pphi t \pphi^{-1}$ with $\taum(t)=t$ so
that $\taum(w)=\tau(\pphi)t\tau(\pphi^{-1}) = \tau(\pphi) \pphi^{-1} w \pphi
\tau(\pphi^{-1})$ (observe that when $U=U(n)$ then $\taum(w)=w^t$ and all of
this becomes clear). Thus $\eta(\alpha(u,v))=\tau(u^{-1}v)$ and
$u^{-1}v=\eta(u,v)$ lie in the same conjugacy class $\calC=\eta(\calDt)$,
so, by lemma \ref{eta-map}, we have indeed $\alpha(u,v)\in\calDt$. From the
remark preceding lemma \ref{D-inv} we see that $Fix(\alpha|_{\calDt})\simeq
\calD \not=\emptyset$.
\end{proof}

On the product $\pdblt$ of $l$ $U_{\Delta}\times U_{\Delta}$-orbits in
$U\times U$, we can therefore define the involution~:

$$\begin{array}{rccl}
\alpha^{(l)} : & \pdblt & \longto & \pdblt \\
& \big((u_1,v_1),\,\ldots\, ,(u_l,v_l)\big) & \longmapsto &
\Big(\big(\tau(v_1),\tau(u_1)\big),\,\ldots\, ,\big(\tau(v_l),\tau(u_l)\big)\Big)
\end{array}$$

\noindent Observe that its fixed point set satisfies $Fix(\alpha^{(l)})\simeq\pdbl$ and is
therefore non-empty. We then have the following result, which says that the
centered Lagrangian problem has a solution if and only if there exists a
solution of the complexified problem which is fixed by $\alpha^{(l)}$~:

\begin{prop}
Let $\calD_1,\,\ldots\, ,\calD_l$ be $l$ $K\times K$-orbits in $U$. For every
$j\in\{1,\,\ldots\, ,l\}$, let $\calC_j$ be the conjugacy class of $\taum(w)w$
where $w$ is any element in $\calD_j$, and let $\calDt_j$ be the
corresponding $U_{\Delta}\times U_{\Delta}$-orbit in $U\times U$ (i.e., such
that $\eta(\calDt_j)=\calC_j$, where $\eta(u,v)=u^{-1}v$). Then there exist $(w_1,\,\ldots\,
,w_l)\in\pdbl$ such that $w_1\ldots w_l=1$ if and only there exist
$((u_1,v_1), \,\ldots\, , (u_l,v_l)) \in \pdblt$ such that $u_1\ldots u_l=1$,
$v_1 \ldots v_l=1$ and $u_j=\tau(v_j)$ for all $j\in\{1,\,\ldots\, ,l\}$ (that
is, $((u_1,v_1), \,\ldots\, , (u_l,v_l)) \in Fix(\alpha^{(l)})$).
\end{prop}

\begin{proof}
For a given $(w_1,\,\ldots\, ,w_l) \in \pdbl ~|~ w_1\ldots w_l=1$, set
$(u_j,v_j):=(\tau(w_j),w_j)$. By lemma \ref{eta-map}, $(u_j,v_j)$ then
belongs to $\calDt_j$ and we have indeed $u_1\ldots u_l=v_1 \ldots v_l = 1$.
Conversely, for $((u_j,v_j))_j \in \pdblt ~|~ u_1\ldots u_l = v_1 \ldots
v_l = 1$ and such that $u_j=\tau(v_j)$ for all $j$, set $w_j:=v_j$. Then
$w_1 \ldots w_l = 1$ and $\taum(w_j)w_j=u_j^{-1}v_j \in \calC_j$, so that,
by corollary \ref{KxK-orbits}, $w_j \in \calD_j$.
\end{proof}

This type of result is exactly why some given problem (A) is called a real
version of another problem (B)~: if $\mathcal{S}_{\C}$ denotes the set of
solutions to problem (B) (we assume that $\mathcal{S}_{\C}\not=\emptyset$)
and $\mathcal{S}_{\R}$ the set of solutions to problem (A), then there
exists an involution $\alpha$ on some space $M\supset \mathcal{S}_{\C}$,
whose fixed point set is non-empty, such that
$\mathcal{S}_{\R}\not=\emptyset$ iff $\mathcal{S}_{\C}\cap
Fix(\alpha)\not=\emptyset$.\\
The question then is~: what is the real version of the unitary problem~?
Given what we have done so far, we see that giving an answer to this
question amounts to defining an involution $\beta^{(l)}$ on $\pconj$ such that
$\beta^{(l)}\circ\eta^{(l)}=\eta^{(l)}\circ\alpha^{(l)}$, where $\eta^{(l)}:\pdblt \to \pconj$
is defined as in proposition \ref{Thompson1}, so that $\eta(Fix(\alpha^{(l)}))
\subset Fix(\beta^{(l)})$, which in particular implies that
$Fix(\beta^{(l)})\not=\emptyset$. The only possibility is then to set, for any
$(w_1,\,\ldots\, ,w_l) \in \pconj$~:
\begin{eqnarray}
\beta^{(l)}(w_1, \,\ldots\, , w_l) \negthickspace & = & \negthickspace
\big(\taum(w_l) \ldots \taum(w_2) \taum(w_1) \tau(w_2) \ldots \tau(w_l),
\,\ldots\, , \,\taum(w_l)\taum(w_{l-1})\tau(w_l),\, \taum(w_l)\big) \quad
\label{def_beta}
\end{eqnarray}
We then have the following result (proposition \ref{Thompson2}), along the
lines of proposition \ref{Thompson1}. As earlier, we see that the group
$K^l$ acts on $Fix(\alpha^{(l)})$ and preserves the relations $u_1 \ldots u_l =
v_1 \ldots v_l =1$. Likewise, $K$ acts diagonally on $Fix(\beta^{(l)})$,
preserving the relation $w_1 \ldots w_l =1$. We may therefore define~:
$$\mathcal{M}_{\calDt}^{\alpha}:=\raisebox{2.5pt}{$\Big\{\big((u_j,v_j)\big)_j\in\pdblt ~|~ u_1\ldots u_l =
v_1\ldots v_l =1\mathrm{~and~} \big((u_j,v_j)\big)_j \in
Fix(\alpha^{(l)})\Big\}$}\bigg/\raisebox{-4pt}{$K^l$}$$
and $$\Mod^{\beta}=\raisebox{2pt}{$\big\{(w_j)_j \in \pconj \,|\, w_1\ldots w_l=1 \mathrm{~and~}
(w_j)_j\in Fix(\beta^{(l)}) \big\}$}\Big/ \raisebox{-3.2pt}{$K$}$$

\noindent We then have~:

\begin{prop}\label{Thompson2}
The map $\eta^{(l)} : \pdblt \to \pconj$ induces a homeomorphism
$\mathcal{M}_{\calDt}^{\alpha}\simeq\Mod^{\beta}$ .
\end{prop}

Again, this is an analog of a result in \cite{AMW}, which justifies that we
may consider our Lagrangian problem a compact version of the Thompson
problem. We may now move on to the main results of this paper.

\subsection{The set of $\s_0$-Lagrangian representations}\label{beta}

Let $\calC_1, \,\ldots\, ,\calC_l$
be $l$ conjugacy classes in $U(n)$ such that there exist $(u_1, \,\ldots\,
, u_l)\in \calC_1\times \cdots \times \calC_l$ satisfying $u_1\ldots u_l=1$.

\begin{defi}
The representation of $\piS$ corresponding to such a $(u_1, \,\ldots\, , u_l)$
is said to be \emph{Lagrangian} if there exist $l$ Lagrangian subspaces
$L_1, \,\ldots\, , L_l$ of $\C^n$ such that, denoting by $\s_j$ the Lagrangian
involution associated to $L_j$, we have $u_j=\s_j\s_{j+1}$ for all
$j\in\{1,\,\ldots\, ,l\}$ (with $\s_{l+1}=\s_1$). It is said to be
$\s_0$\emph{-Lagrangian} if it is Lagrangian with $L_1=L_0:=\R^n\subset\C^n$.
\end{defi}

Recall that two representations $(u_1, \,\ldots\, , u_l)$ and $(v_1, \,\ldots\,
, v_l)$ of $\piS$ are equivalent if and only if there exists a unitary map
$\pphi\in U(n)$ such that $\pphi u_j \pphi^{-1} = v_j$ for all $j\in \{1,
\,\ldots\, , l\}$. Since $\s_{\pphi(L)}=\pphi\s_L\pphi^{-1}$, we have that
any representation equivalent to a Lagrangian one is itself Lagrangian. In
particular, since for any Lagrangian $L\in\lag(n)$ there exists a unitary
map $\pphi\in U(n)$ such that $\pphi(L)=L_0$, we see that a given
representation is Lagrangian if and only if it is equivalent to a
$\s_0$-Lagrangian one.
We now define the map~:

\begin{eqnarray}
\beta :\quad \calC_1 \times \cdots \times \calC_l & \longto & \calC_1 \times
\cdots \times \calC_l \\
(u_1, \,\ldots\, , u_l) & \longmapsto & (\overline{u}_l^{-1} \ldots
\overline{u}_2^{-1} u_1^t \overline{u}_2 \ldots \overline {u}_l, \,\ldots\, ,
\,\overline{u}_l^{-1}u_{l-1}^t\overline{u}_l,\, u_l^t) \label{def_beta_2}
\end{eqnarray}

\noindent (see equation (\ref{def_beta}) in the previous subsection for motivation~: when $U=U(n)$,
$\tau(u)=\overline{u}$). Observe that $\beta$ is an involution (for $l=3$ one
easily sees that $\beta^2=Id$) and that
$Fix(\beta)\not=\emptyset$ (one may for instance pick a diagonal element
$u_j$ in every $\calC_j$ and then $\beta(u_1, \,\ldots\, , u_l)=(u_1, \,\ldots\,
, u_l)$). Also, we have the compatibility relations (see theorem
\ref{lag_locus}) $\beta(\pphi.(u_1, \,\ldots\, , u_l))=\overline{\pphi}.\beta(u_1,
\,\ldots\, , u_l)$ and $\mu\circ\beta(u_1, \,\ldots\, ,
u_l)=\overline{u}_l^{-1}\ldots\overline{u}_1^{-1} = (\overline{\mu(u_1, \,\ldots\, ,
u_l)})^{-1}$, where $\mu$ is the product map $\mu(u_1,\,\ldots\,
u_l)=u_1\ldots u_l$ on $\pconj$. Finally, we consider the Euclidean product on the Lie algebra
$\mathfrak{u}(n)$ given by $(X\,|\,Y)=tr(XY^*)=-tr(XY)$. In particular, the map
$\tau:X\in\mathfrak{u}(n)\mapsto \overline{X}\in \mathfrak{u}(n)$ is an isometry for this
scalar product. We may now state and prove the following characterization of
$\s_0$-Lagrangian representations~:

\begin{thm}\label{s0-Lag}
Given $l$ conjugacy classes $\calC_1, \,\ldots\, , \calC_l$ of unitary
matrices such that there exist $(u_1, \,\ldots\, , $ $u_l)$ $\in \pconj$ satisfying
$u_1 \ldots u_l=1$, the representation of $\piS$ corresponding to such a
$(u_1, \,\ldots\, , u_l)$ is $\s_0$-Lagrangian if and only if $\beta(u_1, \,\ldots\, ,
u_l)=(u_1, \,\ldots\, , u_l)$ (see equation (\ref{def_beta_2}) for a definition of
$\beta$).
\end{thm}

We could as well have defined $\beta$ on $U(n) \times\cdots\times U(n)$ and
obtained a similar result but we deliberately stated our result this way, as
it will be more appropriate to work with the quasi-Hamiltonian space
$\pconj$ in the following.

\begin{proof}[Proof of theorem \ref{s0-Lag}] Let us start with
$(u_1,\, \ldots\, ,u_l)\in Fix(\beta)$, that is~: 
\begin{eqnarray*}
\overline{u}_l^{-1}\ldots\overline{u}_2^{-1}u_1^t\overline{u}_2\dots\overline{u}_l
& = & u_1 \\
\overline{u}_l^{-1}\ldots\overline{u}_3^{-1}u_2^t\overline{u}_3\dots\overline{u}_l & = & u_2 \\
& \vdots & \\
\overline{u}_l^{-1}\ldots\overline{u}_{j+1}^{-1}u_j^t\overline{u}_{j+1}\dots\overline{u}_l
& = & u_j \\
&\vdots & \\
\overline{u}_l^{-1}u_{l-1}^t\overline{u}_l & = & u_{l-1}\\
u_l^t & = & u_l
\end{eqnarray*}
\noindent Then we have $u_l^t=u_l$ (so that $\overline{u}_l=u_l^{-1}$),
$(u_{l-1}u_l)^t=(\overline{u}_l^{-1}u_{l-1}^t\overline{u}_l
u_l)^t=(u_l^tu_{l-1}^t)^t=u_{l-1}u_l$, \ldots, $(u_j\ldots u_l)^t =
(\overline{u}_l^{-1}\ldots\overline{u}_{j+1}^{-1}u_j^t\overline{u}_{j+1}\ldots\overline{u}_l
\ldots\overline{u}_l^{-1}u_{l-1}^t\overline{u}_lu_l)^t = (u_l^tu_{l-1}^t\ldots
u_{j+1}^t u_j^t)^t = u_j\ldots u_l$, \ldots,
and $(u_1\ldots u_l)^t$ $ =
(\overline{u}_l^{-1}\ldots\overline{u}_2^{-1}u_1^t\overline{u}_2\ldots\overline{u}_l
\overline{u}_l^{-1}u_{l-1}^t\overline{u}_lu_l)^t = (u_l^tu_{l-1}^t\ldots
u_2^t u_1^t)^t =
u_1\ldots u_l$. To these $l$ symmetric unitary matrices we can associate, by
proposition \ref{angles}, $l$ Lagrangian subspaces~:
\begin{eqnarray*}
L_1 & := & \{z \in \C^n ~|~ z-(u_1\ldots u_l)\overline{z}=0\} \\
L_2 & := & \{z \in \C^n ~|~ z-(u_2 \ldots u_l)\overline{z}=0\} \\
& \vdots &\\
L_j & := & \{z \in \C^n ~|~ z-(u_j \ldots u_l)\overline{z}=0\} \\
& \vdots & \\
L_{l-1} & := & \{z \in \C^n ~|~ z-(u_{l-1}u_l)\overline{z}=0\}\\
L_l & := & \{z \in \C^n ~|~ z-u_l\overline{z}=0\}
\end{eqnarray*}
\noindent and denote by $\s_j$ the Lagrangian involution associated to
$L_j$. Let us now assume that $(u_1,\,\ldots\, ,u_l)$ satisfy the full hypotheses of
the theorem, that is, that we have $u_1\ldots u_l=1$. Then $L_1=L_0$. Therefore,
by proposition \ref{angles}, since $L_l=\{z-u_l\overline{z}=0\}$, we have
$\s_l\s_0=u_l$, that is, $\s_l\s_1=u_l$. Further, since
$L_2=\{z-(u_2\ldots u_l)\overline{z}=0\}$, we have $\s_2\s_0=u_2\ldots u_l=u_1^{-1}$ hence
$u_1=\s_1\s_2$. Finally, for all $j\in \{2,\,\ldots\, ,l-1\}$, since
$(u_j\ldots u_l)^t=u_j\ldots u_l$, there exists, by
proposition \ref{angles}, a unitary map $\pphi_j \in U(n) ~|~
\pphi_j^t=\pphi_j$ and $\pphi_j^2=u_j\ldots u_l$, and we then have
$\pphi_j(L_0)=L_j$. Set $L'_j=\pphi_2^{-1}(L_j)=L_0$ and
$L'_{j+1}=\pphi_j^{-1}(L_{j+1})$, and denote by $\s'_j$ and $\s'_{j+1}$ the associated
involutions. Then~:
\begin{eqnarray*}
L'_{j+1} & = & \{z ~|~ \pphi_j(z) \in L_{j+1}\} \\
& = & \{z ~|~ \pphi_j(z) - u_{j+1}\ldots u_l\overline{\pphi_j(z)}=0\} \\
& = & \{z ~|~ \pphi_j(z) -
u_{j+1}\ldots u_l\underbrace{\overline{\pphi_j}}_{{=\pphi_j^{-1}}}(\overline{z})=0\} \\
& = & \{z ~|~ z - (\pphi_j^{-1}u_{j+1}\ldots u_l\pphi_j^{-1})\overline{z}=0\}
\end{eqnarray*}
\noindent but $(\pphi_j^{-1}u_{j+1}\ldots u_l\pphi_j^{-1})^t=\pphi_j^{-1}u_{j+1}\pphi_j^{-1}$
since $(\pphi_j^{-1})^t=(\pphi_j^t)^{-1}=\pphi_j^{-1}$ and $(u_{j+1}\ldots
u_l)^t=u_{j+1}\ldots u_l$. Therefore, by proposition \ref{angles}, we have
$\s'_{j+1}\s'_j=\pphi_j^{-1}u_{j+1}\ldots u_l\pphi_j^{-1}$. Since
$\pphi_j^2=u_j\ldots u_l$, we then
have $\pphi_j^{-1}u_{j+1}\ldots u_l\pphi_j^{-1}=\pphi_j^{-1}(u_j^{-1}\pphi_j^2)\pphi_j^{-1}
=\pphi_j^{-1}u_j^{-1}\pphi_j$, therefore
$u_j^{-1}=\pphi_j\s'_{j+1}\s'_j\pphi_j^{-1}=\s_{j+1}\s_j$ since $L_j=\pphi_j(L'_j)$,
$L_{j+1}=\pphi_j(L'_{j+1})$ and $\s_{\pphi(L)}=\pphi\s_L\pphi_{-1}$. Hence
$u_j=\s_j\s_{j+1}$ and the representation of $\pi$ corresponding to
$(u_1,\,\ldots\, ,u_l)$ is $\s_0$-Lagrangian.\\
Conversely, assume that a given representation $(u_1,\,\ldots\, ,u_l)$ is
$\s_0$-Lagrangian. Then $u_l=\s_l\s_0$. Now observe that for any unitary map
$u$, one has $\overline{u}=\s_0 u \s_0$, therefore here
$u_l^t=\overline{u}_l^{-1}=\s_0u_l^{-1}\s_0=\s_0(\s_l\s_0)^{-1}\s_0=\s_0
(\s_0\s_l)\s_0=\s_l\s_0=u_l$. Likewise~:
\begin{eqnarray*}
\overline{u}_l^{-1}u_{l-1}^t\overline{u}_l
& = & (\s_0u_l^{-1}\s_0)(\s_0u_{l-1}^{-1}\s_0)(\s_0u_l\s_0) \\
& = & \s_0(u_l^{-1}u_{l-1}^{-1}u_l)\s_0 \\
& = & \s_0(\s_0\s_l)(\s_l\s_{l-1})(\s_l\s_0)\s_0 \\
& = & \s_{l-1}\s_l \\
& = & u_{l-1}
\end{eqnarray*}
\noindent and so on, until~:
\begin{eqnarray*}
\overline{u}_l^{-1}\ldots\overline{u}_2^{-1}u_1^t\overline{u}_2\dots\overline{u}_l
& = &
\s_0(\s_0\s_l)\ldots(\s_3\s_2)(\s_2\s_1)(\s_2\s_3)\ldots(\s_3\s_0)\s_0 \\
& = & \s_1\s_2 \\
& = & u_1
\end{eqnarray*}
\noindent so that $\beta(u_1,\,\ldots\, ,u_l)=(u_1,\,\ldots\, ,u_l)$.
\end{proof}

We can then characterize those among representations of $\piS$ which are
Lagrangian in the following way~:

\begin{cor}[Characterization of Lagrangian representations]\label{charac}
Suppose that one of the $\calC_j$ is defined by pairwise distinct
eigenvalues and
let $u_1, \,\ldots\, , u_l$ be $l$ unitary matrices such that $u_j\in\calC_j$
and $u_1\ldots u_l=1$. Then there exist $l$ Lagrangian subspaces $L_1,
\,\ldots\, , L_l$ of $\C^n$ such that $u_1=\s_1\s_2, u_2=\s_2\s_3, \,\ldots\, ,
u_l=\s_l\s_1$ (where $\s_j$ is the Lagrangian involution associated to
$L_j$) if and only if $\beta(u_1, \,\ldots\, , u_l)$ is equivalent to $(u_1, \,\ldots\,
, u_l)$ as representations of $\pi$. In this case, if $\psi$ is any unitary
map such that $\beta(u_1, \,\ldots\, , u_l)=(\psi u_1 \psi^{-1}, \,\ldots\, ,
\psi u_l \psi^{-1})$ then $\psi^t=\psi$ and if $\pphi$ is any unitary map
such that $\pphi^t\pphi=\psi$ then the representation of $\pi$ corresponding
to $(\pphi u_1 \pphi^{-1}, \,\ldots\, , \pphi u_l \pphi^{-1})$ is
$\s_0$-Lagrangian.
\end{cor}

\begin{proof} Suppose first that $u_1=\s_1\s_2, \,\ldots\, , u_l=\s_l\s_1$.
Take $\pphi\in U(n) ~|~ \pphi(L_1)=L_0$. Then $\pphi . (u_1, \,\ldots\, ,
u_l)$ is $\s_0$-Lagrangian, hence $\beta(\pphi.(u_1, \,\ldots\, ,
u_l))=\pphi.(u_1, \,\ldots\, , u_l)$, hence $\overline{\pphi}.\beta(u_1,
\,\ldots\, , u_l)=\pphi. (u_1, \,\ldots\, , u_l)$ hence $\beta(u_1, \,\ldots\, ,
u_l)=(\overline{\pphi}^{-1}\pphi). (u_1, \,\ldots\, , u_l) \sim_{U(n)} (u_1,
\,\ldots\, , u_l)$. Observe that $\overline{\pphi}^{-1}\pphi=\pphi^t\pphi$ is
symmetric.\\
Conversely, suppose that $\exists~ \psi \in U(n) ~|~ \beta(u_1, \,\ldots\, ,
u_l)=\psi.(u_1, \,\ldots\, , u_l)$ and assume first that the conjugacy class
$\calC_l$ is defined by pairwise distinct eigenvalues. Write $u_l=vdv^{-1}$
where $d$ is diagonal. Then, since $\beta(u)=\psi.u$, we have in particular
$\psi u_l \psi^{-1} = u_l^t$, from which we obtain $\psi v d v^{-1}
\psi^{-1} = (v^{-1})^t d^t v^t = (v^t)^{-1} d v^t$, so that $(v^t \psi v) d
(v^t \psi v)^{-1} = d$. Since $d$ is diagonal with pairwise distinct
elements, $v^t \psi v$ is itself diagonal and therefore symmetric, so that
$\psi$ is symmetric. If now it is a different $\calC_j$ which is defined by
pairwise distinct eigenvalues, say $\calC_{l-1}$, then consider the
representation $(u_l,u_1, \,{\ldots}\, ,u_{l-1})$~: it is indeed a
representation of $\pi$ since the relation $u_1 \ldots u_l =1$ is invariant
by circular permutation (as can be seen by conjugating by $u_l$) and via
this transformation $\psi.(u_1, \,\ldots\, , u_l)$ is sent to
$\psi.(u_l,u_1, \,\ldots\, ,u_{l-1})$. The representation $(u_l,u_1,
\,\ldots\, , u_{l-1})$ is Lagrangian iff $(u_1, \,\ldots\, , u_l)$ is
Lagrangian. We can define a corresponding $\beta$ accordingly and proceed as
above to show that $\psi$ is indeed symmetric. \\
Now, to conclude, let $\pphi$ be any unitary map such that $\pphi^t\pphi=\psi$ (such a map
always exists by proposition \ref{angles}). Starting from $\beta(u_1,
\,\ldots\, , u_l)=\psi.(u_1, \,\ldots\, , u_l)$, we obtain
$(\pphi^t)^{-1}.\beta(u)=\pphi.u$, hence $\beta(\pphi.u)=\pphi.u$ so that, by
theorem \ref{s0-Lag}, $\pphi.(u_1, \,\ldots\, , u_l)$ is $\s_0$-Lagrangian.
Hence $(u_1, \,\ldots\, , u_l)$ is Lagrangian, with $L_1=\pphi^{-1}(L_0)$.
\end{proof}

Before passing on to studying Lagrangian representations in the moduli
space, we would like to point out that if a representation $u$ is
irreducible then so is $\beta(u)$ and, more interestingly maybe, that it is
possible to characterize Lagrangian representations with arbitrarily fixed
first Lagrangian $L_1$ in a way similar to theorem \ref{s0-Lag}. In order to
do so, we define, for a given Lagrangian subspace $L_1$,
the involution $\beta_{L_1}(u_1, u_2 , u_3):= (\s_1
u_3^{-1} u_2^{-1} u_1^{-1} u_2 u_3 \s_1, \s_1 u_3^{-1} u_2^{-1} u_3 \s_1,
\s_1 u_3^{-1} \s_1)$ (remember that when $L_1=L_0$, $\s_0 u \s_0 =
\overline{u}$). If we write $L_1=\pphi(L_0)$ for some $\pphi \in U(n)$, we
obtain $\beta_{L_1}(u)=(\pphi\pphi^t).\beta(u)$ (this does not depend on the
choice of $\pphi$ such that $\pphi(L_0)=L_1$ as seen from the argument used
in proposition \ref{angles}). Finally, it was proved in \cite{FMS} that when $n=2$ and
$l=3$, every (two-dimensional) unitary representation of
$\pi_1(S^2\bs\{s_1,s_2,s_3\})$ is Lagrangian~: this is because in this case
the moduli space is a single point (it is zero-dimensional and connected),
so that the submanifold consisting of Lagrangian representations is the
point itself (see \cite{FW} for dimensions of moduli spaces of
representations). As a matter of fact, we believe that the
characterization of Lagrangian representations as representations $u$
satisfying $\beta(u) \sim_{U(n)} u$ is true even without the (generic)
assumption made on the
$\calC_j$ but we have been unable to prove it so far. One would only need to
show that if $\beta(u)=\psi.u$ for some $\psi \in U(n)$ then there exists
such a $\psi$ which is symmetric. In the remainder of this paper, we will assume
that one of the $\calC_j$ is defined by pairwise distinct eigenvalues, so
that corollary \ref{charac} holds.

\begin{rk}[Addendum - 26.07.06]
As a matter of fact, corollary \ref{charac} does hold without any assumption on the conjugacy classes $\calC_j$ and a proof of this is availablle in \cite{mythesis}. 
\end{rk}

\subsection{Lagrangian representations in the moduli
space}\label{lag_submanifold}

Recall from section \ref{moduli} that the moduli space of unitary
representations of $\pi=\piS$ is the quasi-Hamiltonian quotient
$\Mod=\mu^{-1}(\{1\})/U(n)$ where $\mu:\pconj\to U(n)$ is the product map.
Since the involution $\beta$ we constructed on $\pconj$ in \ref{beta}
satisfies $\beta\circ\mu = \taum\circ\mu$ (where $\tau(u)=\overline{u}$ on
$U(n)$), $\beta$ preserves $\mu^{-1}(\{1\})$ and since
$\beta(\pphi.u)=\tau(\pphi).\beta(u)$, $\beta$ induces an involution
$\hat{\beta}$ on $\Mod=\mu^{-1}(\{1\})/U(n)$ given by
$\hat{\beta}([u])=[\beta(u)]$. Observe that if $\beta^{(l)}$ is defined as in
the end of the previous subsection by $\beta_{L}=(\pphi\pphi^t).\beta$
(where $\pphi \in U(n)$ satisfies $\pphi(L_0)=L$) then
$\widehat{\beta^{(l)}}=\hat{\beta}$. Furthermore, if $[u]\in\Mod$ is the
equivalence class of a unitary representation of $\pi$, then it is
Lagrangian if and only if any of its representatives is Lagrangian (for, if
$u_j=\s_j\s_{j+1}$, then $\pphi u_j \pphi^{-1} = \s'_j\s'_{j+1}$ where
$L'_j=\pphi(L_j)$, for any $\pphi\in U(n)$). Corollary \ref{charac} then
shows that a given $[u]\in\Mod$ is Lagrangian if and only if
$\hat{\beta}([u]) = [u]$. We then have the following result, which is a
direct consequence of theorem \ref{lag_locus}.

\begin{thm}\label{lag_moduli}
The set of equivalence classes of Lagrangian representations of $\pi=\piS$
is exactly $Fix(\hat{\beta})$. It is a Lagrangian submanifold of the moduli
space $\Mod=\Hom_{\calC}(\pi,U(n))/U(n)$ of unitary representations of $\pi$
(in particular it is always non-empty).
\end{thm}

To apply theorem \ref{lag_locus}, the only condition left to check is that
$\beta^*\w=-\w$, where $\w$ is the $2$-form defining the quasi-Hamiltonian
structure on $\pconj$ described in section \ref{moduli}. Actually, we also
need to check that $Fix(\hat{\beta})\not=\emptyset$. As indicated before theorem
\ref{convexity},
this is always true for an involution $\beta$ which satisfies
the hypotheses of theorem \ref{lag_locus} and which has fixed points itself,
but since this paper does not
contain a proof of this fact, we instead refer to theorem 1 of \cite{FW},
which we state here.

\begin{thm}\emph{\cite{FW}}
Let $\calC_1, \,\ldots\, , \calC_l$ be $l\geq 1$ conjugacy classes in $U(n)$
such that there exist $(u_1, \,\ldots\, ,$ $u_l)$ $\in \pconj$ satisfying
$u_1\ldots u_l =1$. Then there exist $l$ Lagrangian subspaces $L_1,
\,\ldots\, , L_l$ of $\C^n$ such that $\s_j\s_{j+1}\in\calC_j$ for all $j\in
\{1, \,\ldots\, , l\}$, where $\s_j$ is the Lagrangian involution associated
with $L_j$ and where $\s_{l+1}=\s_1$.
\end{thm}

\noindent This shows that $Fix(\beta)\cap\mu^{-1}(\{1\})\not=\emptyset$ as
one can construct, from the Lagrangian representation $(\s_j\s_{j+1})_j$ whose existence is
guaranteed by the theorem, a $\s_0$-Lagrangian representation
$(\s'_j\s'_{j+1})_j\in Fix(\beta)$ by applying $\pphi\in U(n)$ such that
$\pphi(L_1)=L_0$.

\begin{proof}[Proof of theorem \ref{lag_moduli}]

As observed, we only have to check that $\beta^* \w = - \w$. We prove it by
induction on $l$. For $l=1$, we have, for any $X,Y \in \mathfrak{u}$
(denoting $[X]_u = X.u - u.X \in T_u\calC_1$), $$\w_u([X]_u,[Y]_u) =
\frac{1}{2} \big((Ad u.X \,|\, Y) - (Ad u.Y \,|\, X)\big)$$ as well as $\beta(u)=\tau(u^{-1})$
and $T_u\beta.[X]_u = [\tau(X)]_{\tau(u^{-1})}$. Therefore~:

\begin{eqnarray*}
(\beta^*\w)_u\big([X]_u,[Y]_u\big) & = & 
\w_{\beta(u)}\big(T_u\beta.[X]_u,T_u\beta.[Y]_u\big) \\
& = & \frac{1}{2}\, \Big(\big(Ad\,\tau(u^{-1}).\tau(X)\,|\,\tau(Y)\big) -
\big(Ad\,\tau(u^{-1}).\tau(Y)\,|\,\tau(X)\big)\Big) \\
& = & \frac{1}{2}\, \Big(\big(\tau(Ad\,u^{-1}.X)\,|\,\tau(Y)\big) -
\big(\tau(Ad\,u^{-1}.Y)\,|\,\tau(X)\big)\Big)
\end{eqnarray*}

\noindent Since $\tau$ is an isometry for $(.\,|\,.)$, we then have~:

\begin{eqnarray*}
(\beta^*\w)_u\big([X]_u,[Y]_u\big)
& = & \frac{1}{2}\, \big((Ad\,u^{-1}.X\,|\,Y) -
(Ad\,u^{-1}.Y\,|\,X)\big)\\
& = & \frac{1}{2}\, \big((X\,|\,Ad\,u.Y) - (Y\,|\,Ad\,u.X)\big) \\
& = & -\w_u\big([X]_u,[Y]_u\big)
\end{eqnarray*}

\noindent To complete the induction, we will use the following lemma, which
is general in nature and can be
used to construct form-reversing involutions on quasi-Hamiltonian spaces.

\begin{lem}\label{anti-inv}
Let $(M_1,\w_1,\mu_1: M_1 \to U)$ and $(M_2,\w_2,\mu_2: M_2 \to U)$ be two
quasi-Hamiltonian $U$-spaces. Let $\tau$ be an involutive automorphism of
$(U,(. \,|\, .))$ and let $\beta_i$ be an involution on $M_i$ satisfying~:

\begin{enumerate}
\item[(i)] $\beta_i^* \w_i = -\w_i$
\item[(ii)] $\beta_i(u.x_i) = \tau(u).\beta_i(x_i)$ for all $u \in U$ and all
$x_i \in M_i$
\item[(iii)] $\mu_i \circ \beta_i = \taum \circ \mu_i$
\end{enumerate}

\noindent Consider the quasi-Hamiltonian $U$-space $(M:= M_1 \times M_2, \w:= \w_1
\oplus \w_2 + (\mu_1^*\theta^L \wedge \mu_2^* \theta^R), \mu:=\mu_1 \cdot
\mu_2)$ (with respect to the diagonal action of $U$) and the map~:

$$\begin{array}{rccl}
\beta:= \big( (\mu_2 \circ \beta_2).\beta_1, \beta_2 \big): & M & \longto & M \\
& (x_1,x_2) & \longmapsto & \big( (\mu_2\circ\beta_2(x_2)).\beta_1(x_1),
\beta_2(x_2)\big)
\end{array}$$

\noindent Then $\beta$ is an involution on $M$ satisfying~:

\begin{enumerate}
\item[(i)] $\beta^* \w = -\w$
\item[(ii)] $\beta(u.x) = \tau(u).\beta(x)$ for all $u \in U$ and all
$x \in M$
\item[(iii)] $\mu \circ \beta = \taum \circ \mu$
\end{enumerate}

\end{lem}

\noindent We postpone the proof of the lemma and give the end of the proof of theorem
\ref{lag_moduli}. To complete the induction, all one has to do is check that our involution
$\beta=\beta^{(l)}$ (see (\ref{def_beta})) on the product $\pconj$ of $l$ conjugacy classes is indeed obtained like
in the lemma starting from the form-reversing involution
$\beta^{(1)}:=\taum: u \to u^t$ on each single
conjugacy class. This is easily checked since on $\calC_1\times\calC_2$~:

\begin{eqnarray*}
\beta^{(2)}(u_1,u_2) & = & (\overline{u_2}^{-1} \, u_1^t \,
\overline{u_2}\, , u_2^t) \\
& = & (u_2^t . u_1^t \, , u_2^t) \\
& = & \big( (\mu_2\circ\beta^{(1)}(u_2)).\beta^{(1)}(u_1), \beta^{(1)}(u_2) \big)
\end{eqnarray*}

\noindent and on $\calC_1\times(\calC_2\times\calC_3)$~:

\begin{eqnarray*}
\beta^{(3)}(u_1,u_2,u_3) & = & (\overline{u_3}^{-1} \, \overline{u_2}^{-1} \,
u_1^t \, \overline{u_2} \, \overline{u_3}\, , \overline{u_3}^{-1}\, u_2^t
\, \overline{u_3}\, , u_3^t) \\
& = & \big( (u_2 u_3)^t . u_1^t \, , u_3^t.u_2^t\, , u_3^t \big) \\
& = & \Big( \big(
(\mu_2\cdot\mu_3)\circ\beta^{(2)}(u_2,u_3)\big).\beta^{(1)}(u_1),
\beta^{(2)}(u_2,u_3) \Big)
\end{eqnarray*}

\noindent and so on. It is of course the very form of the involution $\beta$
which inspired the formulation of the lemma.
\end{proof}

\begin{proof}[Proof of lemma \ref{anti-inv}]
First, we have~:
\begin{eqnarray*}
\beta(\beta(x_1,x_2)) & = & \bigg( \Big( \mu_2 \circ \beta_2
\big(\beta_2(x_2)\big)
\bigg).\beta_1\Big(\big(\mu_2\circ\beta_2(x_2)\big).\beta_1(x_1)\Big),
\beta_2\big(\beta_2(x_2))\bigg) \\
& = & \bigg(
\big(\mu_2(x_2)\big).\Big(\tau\big(\underset{=\taum\circ\mu_2}{\underbrace{\mu_2\circ\beta_2}}(x_2)\big).
\beta_1\big(\beta_1(x_1)\big)\Big), x_2\bigg) \\
& = & \Big(\big(\mu_2(x_2)\big)\big(\mu_2(x_2)\big)^{-1}.x_1,x_2\Big) \\
& = & (x_1,x_2)
\end{eqnarray*}

\noindent so that $\beta$ is indeed an involution. Second~:

\begin{eqnarray*}
\beta(u.x_1,u.x_2) & = & \big( \mu_2\circ\beta_2(u.x_2).\beta_1(u.x_1),
\beta_2(u.x_2)\big) \\
& = & \Big(
\underset{=\tau(u)\mu_2\big(\beta_2(x_2)\big)\tau(u)^{-1}}{\underbrace{\mu_2\big(\tau(u).\beta_2(x_2)\big)}}
.\,\big(\tau(u).\beta_1(x_1)\big), \tau(u)\beta_2(x_2)\Big) \\
& = & \bigg( \tau(u).\Big(\big(\mu_2\circ\beta_2(x_2)\big).\beta_1(x_1)\Big),
\tau(u).\beta_2(x_2)\bigg) \\
& = & \tau(u).\beta(x_1,x_2)
\end{eqnarray*}

\noindent and~:

\begin{eqnarray*}
\mu \circ \beta(x_1,x_2) & = &
\mu_1\Big(\big(\mu_2\circ\beta_2(x_2)\big).\beta_1(x_1)\Big)
\mu_2\big(\beta_2(x_2)\big) \\
& = &
\Big(\mu_2\circ\beta_2(x_2)\mu_1\circ\beta_1(x_1)
\big(\mu_2\circ\beta_2(x_2)\big)^{-1}\Big) \Big(\mu_2\circ\beta_2(x_2)\Big) \\
& = & \taum\circ\mu_2(x_2) \taum\circ\mu_1(x_1) \\
& = & \taum\circ(\mu_1\cdot\mu_2)(x_1,x_2) \\
& = & \taum\circ\mu(x_1,x_2)
\end{eqnarray*}

\noindent So the only thing left to prove is that $\beta^*\w =-\w$. Let us
start by computing $T\beta$. For all $(x_1,x_2)\in M$, and all
$(v_1,v_2):=\frac{d}{dt}|_{t=0} (x_1(t), x_2(t))$ (where $x_i(0)=x_i$), one has~:

\begin{eqnarray*}
T_{(x_1,x_2)}\beta.(v_1,v_2) & = &
\frac{d}{dt}|_{t=0}\Big(\big(\mu_2\circ\beta_2(x_2)\big) .
\beta_1\big(x_1(t)\big), \beta_2\big(x_2(t)\big)\Big) \qquad \qquad \qquad
\qquad \qquad \qquad  \qquad \qquad 
\end{eqnarray*}
\begin{eqnarray*}
= & \bigg(
\big(\mu_2\circ\beta_2(x_2)\big).\left\{\Big(\theta^L_{\mu_2\circ\beta_2(x_2)}
\big(T_{x_2}(\mu_2\circ\beta_2).v_2\big)\Big)^{\sharp}_{\beta_1(x_1)} +
T_{x_1}\beta_1.v_1\right\}, T_{x_2}\beta_2.v_2\bigg) &
\end{eqnarray*}

\noindent Recall indeed that if a Lie group $U$ acts on a manifold $M$ then~:
$$ \frac{d}{dt}|_{t=0}\ (u_t.x_t) = u_0.X^{\sharp}_{x_0} +
u_0.\Big(\frac{d}{dt}|_{t=0} \ x_t\Big) $$
\noindent where $X \in \mathfrak{u}=Lie(U)$ is such that $u_t=u_0\exp(tX)$
for all $t$, that is~: $$X=u_0^{-1}.\Big(\frac{d}{dt}|_{t=0} \ u_t\Big) =
\theta^L_{u_0}\Big(\frac{d}{dt}|_{t=0} \ u_t\Big)$$
\noindent Let us now
compute $\beta^*(\w_1\oplus\w_2)$. We obtain, for all $(x_1,x_2) \in M$ and
all $(v_1,v_2),(w_1,w_2) \in T_{(x_1,x_2)}M$~:\vskip 15pt
$\big(\beta^*(\w_1\oplus\w_2)\big)_{(x_1,x_2)}\big(
(v_1,v_2),(w_1,w_2) \big)$
\begin{eqnarray}
& = & (\w_1)_{\big(\mu_2\circ\beta_2(x_2)\big).\beta_1(x_1)} \Bigg(
\big(\mu_2\circ\beta_2(x_2)\big).\left\{\Big(\theta^L_{\mu_2\circ\beta_2(x_2)}
\big(T_{x_2}(\mu_2\circ\beta_2).v_2\big)\Big)^{\sharp}_{\beta_1(x_1)} +
\underset{(A)}{\underbrace{T_{x_1}\beta_1.v_1}}
\right\}, \label{firstline} \\
& & \qquad\qquad\qquad\qquad\qquad\big(\mu_2\circ\beta_2(x_2)\big).\left\{\Big(\theta^L_{\mu_2\circ\beta_2(x_2)}
\big(T_{x_2}(\mu_2\circ\beta_2).w_2\big)\Big)^{\sharp}_{\beta_1(x_1)} +
\underset{(A)}{\underbrace{T_{x_1}\beta_1.w_1}}\right\}\Bigg)
\label{secondline}
\end{eqnarray}
\begin{eqnarray*}
+ \ \underset{(B)}{\underbrace{(\w_2)_{\beta_2(x_2)} (T_{x_2}\beta_2.v_2,
T_{x_2}\beta_2.w_2)}}
\qquad\qquad\qquad\qquad\qquad\qquad\qquad\qquad\qquad\qquad\quad
\end{eqnarray*}

\noindent Since $\w_1$ is $U$-invariant, we can drop the terms
$\mu_2\circ\beta_2(x_2) \in U$ appearing on lines (\ref{firstline}) and
(\ref{secondline}). Further, since $\beta^*\w_1=-\w_1$ and
$\beta^*\w_2=-\w_2$, we have, by the $l=1$ case~:
\begin{eqnarray}
(A) \, + \, (B) & = & -(w_1)_{x_1}(v_1,w_1) - (w_2)_{x_2}(v_2,w_2) \label
{A+B}\\
& = & - (\w_1 \oplus \w_2)_{(x_1,x_2)}\big( (v_1,v_2),(w_1,w_2)\big)
\end{eqnarray}

\noindent The remaining terms on lines (\ref{firstline}) and
(\ref{secondline}) then are~:
\begin{eqnarray}
& & (w_1)_{\beta_1(x_1)} \bigg(\Big(\theta^L_{\mu_2\circ\beta_2(x_2)} \big(
T_{x_2}(\mu_2\circ\beta_2).v_2\big)\Big)^{\sharp}_{\beta_1(x_1)} \, ,
T_{x_1}\beta_1 .w_1\bigg) \label{firstagain}\\
& + & \ (w_1)_{\beta_1(x_1)} \bigg( 
T_{x_1}\beta_1.v_1\, , \Big(\theta^L_{\mu_2\circ\beta_2(x_2)} \big(
T_{x_2}(\mu_2\circ\beta_2).w_2\big)\Big)^{\sharp}_{\beta_1(x_1)}\bigg)
\label{secondagain}\\
& + & \ (w_1)_{\beta_1(x_1)} \bigg(\Big(\theta^L_{\mu_2\circ\beta_2(x_2)} \big(
T_{x_2}(\mu_2\circ\beta_2).v_2\big)\Big)^{\sharp}_{\beta_1(x_1)} \, ,
\Big(\theta^L_{\mu_2\circ\beta_2(x_2)} \big(
T_{x_2}(\mu_2\circ\beta_2).w_2\big)\Big)^{\sharp}_{\beta_1(x_1)}
\bigg) \label{thirdagain} 
\end{eqnarray}

\noindent and we notice that each of these three terms is of the form
$\iota_{X^{\sharp}}\w_1 = \frac{1}{2}\mu_1^* (\theta^L + \theta^R \,|\, X)$
for some $X \in \mathfrak{u}$. To facilitate the computations, we set, for
$i=1,2$~:

\begin{eqnarray*}
g_i & := & \mu_i \circ \beta_i (x_i) \in U \\
\zeta_i & := & T_{x_i} (\mu_i\circ\beta_i) . v_i \in
T_{\mu_i\circ\beta_i(x_i)}U=T_{g_i} U \\
\eta_i & := & T_{x_i} (\mu_i\circ\beta_i) . w_i \in
T_{\mu_i\circ\beta_i(x_i)}U=T_{g_i} U \\
\end{eqnarray*}

\noindent We can then rewrite lines (\ref{firstagain}), (\ref{secondagain})
and (\ref{thirdagain}) under the form~:
\begin{eqnarray}
& & \frac{1}{2} \big( \underset{(1)}{\underbrace{\theta^L_{g_1}(\eta_1)}} +
\underset{C}{\underbrace{\theta^R_{g_1}(\eta_1)}} \,|\,
\theta^L_{g_2}(\zeta_2) \big) \label{one}\\
& & - \frac{1}{2} \big( \underset{(2)}{\underbrace{\theta^L_{g_1}(\zeta_1)}} +
\underset{D}{\underbrace{\theta^R_{g_1}(\zeta_1)}} \,|\,
\theta^L_{g_2}(\eta_2) \big) \label{two}\\
& & +\frac{1}{2} \Big(
\theta^L_{g_1}\big(\theta^L_{g_2}(\eta_2).g_1-g_1.\theta^L_{g_2}(\eta_2)\big)
+ \theta^R_{g_1}\big(\theta^L_{g_2}(\eta_2).g_1-g_1.\theta^L_{g_2}(\eta_2)\big)
\,|\, \theta^L_{g_2}(\zeta_2)\Big) \label{three}
\end{eqnarray}

\noindent where the expression for the last term follows from the
equivariance of $\mu_1$~:
$$ T_{\beta_1(x_1)}\mu_1 .
\big(\theta^L_{\mu_2\circ\beta_2(x_2)}(T_{x_2}(\mu_2\circ\beta_2).w_2)\big)^{\sharp}_{\beta_1(x_1)}
=
\big(\theta^L_{\mu_2\circ\beta_2(x_2)}(T_{x_2}(\mu_2\circ\beta_2).w_2)\big)^{\sim}_{\mu_1\circ\beta_1(x_1)}
= \big(\theta^L_{g_2}(\eta_2)\big)^{\sim}_{g_2}
$$

\noindent (where $X^{\sim}_u = X.u - u.X$ is the value at $u$ of the
fundamental vector field associated to $X\in \mathfrak{u}$ by the action of
$U$ on itself by conjugation). We can simplify the expression in (\ref{three})
further by using the definition of $\theta^L$ and $\theta^R$ and the
$Ad$-invariance of $(. \,|\, .)$~:

\begin{eqnarray}
(\ref{three}) & = & \frac{1}{2} \big( Ad\, g_1^{-1}.\theta^L_{g_2}(\eta_2) -
Ad\, g_1.\theta^L_{g_2}(\eta_2) \,|\, \theta^L_{g_2}(\zeta_2) \big) \\
& = & \frac{1}{2} \big( \theta^L_{g_2}(\eta_2) \,|\, Ad\, g_1 . \theta^L_{g_2}(\zeta_2) \big)
- \frac{1}{2}\big(Ad\, g_1.\theta^L_{g_2}(\eta_2) \,|\,
\theta^L_{g_2}(\zeta_2) \big) \label{threeagain}
\end{eqnarray}

\noindent Let us now compute $\beta^*(\mu_1^*\theta^L \wedge \mu_2^*\theta^R)$.
\begin{eqnarray}
& & \big(\beta^*(\mu_1^*\theta^L \wedge
\mu_2^*\theta^R)\big)_{(x_1,x_2)}\big( (v_1,v_2),(w_1,w_2) \big) \\
& = & (\mu_1^*\theta^L \wedge
\mu_2^*\theta^R)_{\big(\mu_2\circ\beta_2(x_2).\beta_1(x_1),\beta_2(x_2)\big)}
\big(T_{(x_1,x_2)}\beta . (v_1,v_2) , T_{(x_1,x_2)}\beta . (w_1,w_2)\big) \\
& = & \frac{1}{2} \bigg(
\theta^L_{\underset{g_2\mu_1(\beta_1(x_1))g_2^{-1}}{\underbrace{\mu_1(g_2.\beta_1(x_1)}}}
\Big( T_{g_2.\beta_1(x_1)}\mu_1.\big(\mu_2\circ\beta_2(x_2)\big).\left\{
\big(\theta^L_{g_2}(\zeta_2)\big)^{\sharp}_{\beta_1(x_1)} +
T_{x_1}\beta_1.v_1 \right\}\Big) \,|\, \theta^R_{g_2}(\eta_2) \bigg)
\label{above}
\end{eqnarray}
\vskip -10pt
\begin{eqnarray*}
& & -\frac{1}{2} 
\bigg(\theta^L_{\mu_1(g_2.\beta_1(x_1)}
\Big( T_{g_2.\beta_1(x_1)}\mu_1.\big(\mu_2\circ\beta_2(x_2)\big).\left\{
\big(\theta^L_{g_2}(\eta_2)\big)^{\sharp}_{\beta_1(x_1)} +
T_{x_1}\beta_1.w_1 \right\}\Big) \,|\, \theta^R_{g_2}(\zeta_2) \bigg)
\end{eqnarray*}

\noindent Since $\mu_1$ is equivariant, we have, for any $v \in
T_{\beta_1(x_1)} M_1$~:
$$T_{g_2.\beta_1(x_1)}\mu_1.\big(\mu_2\circ\beta_2(x_2)\big).v =
\big(\mu_2\circ\beta_2(x_2)\big).\Big(T_{\beta_1(x_1)}\mu_1.v\Big)$$
\noindent where the action in the right side term is conjugation. We then
have~:
\begin{eqnarray}
(\ref{above}) & = & \frac{1}{2}\bigg( \theta^L_{g_2g_1g_2^{-1}}\Big(
g_2. \Big(T_{\beta_1(x_1)}\mu_1 .
\Big(\big(\theta^L_{g_2}(\zeta_2)\big)^{\sharp}_{\beta_1(x_1)} +
T_{x_1}\beta_1.v_1 \Big) \Big).g_2^{-1} \Big) \,|\, \theta^R_{g_2}(\eta_2) \bigg)
\end{eqnarray}
\begin{eqnarray*}
& & \qquad -\frac{1}{2}\bigg( \theta^L_{g_2g_1g_2^{-1}}\Big(
g_2. \Big(T_{\beta_1(x_1)}\mu_1 .
\Big(\big(\theta^L_{g_2}(\eta_2)\big)^{\sharp}_{\beta_1(x_1)} +
T_{x_1}\beta_1.w_1 \Big) \Big).g_2^{-1} \Big) \,|\, \theta^R_{g_2}(\zeta_2) \bigg)
\end{eqnarray*}
\begin{eqnarray}
& = & \frac{1}{2} \Big( g_2 g_1^{-1} g_2^{-1} g_2 .
\big(\theta^L_{g_2}(\zeta_2).g_1 - g_1.\theta^L_{g_2}(\zeta_2)\big) .
g_2^{-1} \,|\, \theta^R_{g_2}(\eta_2) \Big) + \frac{1}{2} \Big( g_2
\underset{=\theta^L_{g_1}(\zeta_1)}{\underbrace{g_1^{-1}
g_2^{-1} g_2 . \zeta_1}} . g_2^{-1} \,|\, \theta^R_{g_2}(\eta_2) \Big)
\end{eqnarray} \vskip -10pt
\begin{eqnarray*}
& & - \frac{1}{2} \Big( g_2 g_1^{-1} g_2^{-1} g_2 .
\big(\theta^L_{g_2}(\eta_2).g_1 - g_1.\theta^L_{g_2}(\eta_2)\big) .
g_2^{-1} \,|\, \theta^R_{g_2}(\zeta_2) \Big) + \frac{1}{2} \Big( g_2
\underset{=\theta^L_{g_1}(\eta_1)}{\underbrace{g_1^{-1}
g_2^{-1} g_2 . \eta_1}} . g_2^{-1} \,|\, \theta^R_{g_2}(\zeta_2) \Big)
\end{eqnarray*}
\begin{eqnarray}
& = & \frac{1}{2} \Big( Ad\, g_2 Ad\, g_1^{-1} .
\theta^L_{g_2}(\zeta_2) \,|\, \theta^R_{g_2}(\eta_2) \Big) - \frac{1}{2}
\Big( Ad\, g_2 . \theta^L_{g_2}(\zeta_2) \,|\, \theta^R_{g_2}(\eta_2) \Big)
+ \frac{1}{2} \Big( Ad\, g_2 . \theta^L_{g_1}(\zeta_1) \,|\, \theta^R_{g_2}(\eta_2) \Big)
\end{eqnarray}
\begin{eqnarray*}
& & - \frac{1}{2} \Big( Ad\, g_2 Ad\, g_1^{-1} .
\theta^L_{g_2}(\eta_2) \,|\, \theta^R_{g_2}(\zeta_2) \Big) + \frac{1}{2}
\Big( Ad\, g_2 . \theta^L_{g_2}(\eta_2) \,|\, \theta^R_{g_2}(\zeta_2) \Big)
- \frac{1}{2} \Big( Ad\, g_2 . \theta^L_{g_1}(\eta_1) \,|\, \theta^R_{g_2}(\zeta_2) \Big)
\end{eqnarray*}
\begin{eqnarray}
& = & \underset{(3')}{\underbrace{\frac{1}{2} \Big( \theta^L_{g_2}(\zeta_2) \,|\, Ad\, g_1 .
\theta^L_{g_2}(\eta_2)\Big)}} -\underset{(4)}{\underbrace{\frac{1}{2}\Big( \theta^L_{g_2}(\zeta_2) \,|\,
\theta^L_{g_2}(\eta_2)}}\Big) + \underset{(2')}{\underbrace{\frac{1}{2} \Big( \theta^L_{g_1}(\zeta_1) \,|\,
\theta^L_{g_2}(\eta_2)\Big)}} \qquad\qquad\qquad\qquad \label{onetwothreeprime}
\end{eqnarray}\vskip -10pt
\begin{eqnarray*}
&  & - \underset{(3')}{\underbrace{\frac{1}{2} \Big( \theta^L_{g_2}(\eta_2) \,|\, Ad\, g_1 .
\theta^L_{g_2}(\zeta_2)\Big)}} +\underset{(4')}{\underbrace{\frac{1}{2}\Big( \theta^L_{g_2}(\eta_2) \,|\,
\theta^L_{g_2}(\zeta_2)}}\Big) - \underset{(1')}{\underbrace{\frac{1}{2} \Big( \theta^L_{g_1}(\eta_1) \,|\,
\theta^L_{g_2}(\zeta_2)\Big)}} \qquad\qquad\qquad\qquad 
\end{eqnarray*}

\noindent (to obtain this last expression, one uses the $Ad$-invariance of
$(.\,|\, .)$ and the fact that $Ad\, g_2^{-1} \circ \theta^R_{g_2} =
\theta^L_{g_2}$). Observe that $(4)$ and $(4')$ cancel in the above
expression. Likewise, $(1')$, $(2')$ and $(3')$ in (\ref{onetwothreeprime})
cancel respectively with $(1)$, $(2)$ in (\ref{one}) and (\ref{two}) and with
(\ref{threeagain}) when computing the sum
$\beta^*(\w_1\oplus\w_2) + \beta^*(\mu_1^*\theta^L \wedge \mu_2^*\theta^R)$.
The non-vanishing terms in this sum are therefore $(A)$ and $(B)$ from
(\ref{A+B}) and $(C)$ and $(D)$ from (\ref{one}) and (\ref{two}), so that~:
\begin{eqnarray}
& & (\beta^*\w)_{x}(v,w) \\
& = & \big( \beta^*(\w_1\oplus\w_2)\big)_x (v,w) +
\big(\beta^*(\mu_1^*\theta^L \wedge \mu_2^*\theta^R)\big)_x (v,w) \\
& = & (A) + (B) + (C) + (D) \\
& = & -(\w_1\oplus\w_2)_x (v,w) - \frac{1}{2}\Big(
\big(\theta^R_{g_1}(\zeta_1) \,|\, \theta^L_{g_2}(\eta_2)\big)
- \big(\theta^R_{g_1}(\eta_1) \,|\, \theta^L_{g_2}(\zeta_2)\big) \Big)
\label{prefinal}
\end{eqnarray}

\noindent But $\mu_i\circ\beta_i = \taum\circ\mu_i$, so that~:
\begin{eqnarray}
\big(\theta^R_{g_1}(\zeta_1) \,|\, \theta^L_{g_2}(\eta_2)\big) & = &
\big(\theta^R_{\mu_1\circ\beta_1(x_1)}(T_{x_1}(\mu_1\circ\beta_1).v_1) \,|\,
\theta^L_{\mu_2\circ\beta_2(x_2)}(T_{x_2}(\mu_2\circ\beta_2).w_2)\big) \\
& = & \big( \theta^R_{\taum\circ\mu_1(x_1)} (T_{x_1}(\taum\circ\mu_1).v_1)
\,|\, \theta^L_{\taum\circ\mu_2(x_2)} (T_{x_2}(\taum\circ\mu_2).w_2)
\big) \label{step}
\end{eqnarray}

\noindent and $\taum=Inv\circ\tau$, where $Inv: u\mapsto u^{-1}$ is
inversion on $U$, so $T_u\taum . \xi = - \taum(u).(T_u \tau.\xi).\taum(u)$.
Hence~:

\begin{eqnarray*}
\theta^R_{\taum(u)}(T_u \taum .\xi) & = & \theta^R_{\taum(u)}
\big(-\taum(u).(T_u\tau.\xi).\taum(u)\big)\\
& = & -\taum(u).(T_u \tau .\xi) \\
& = & -\theta^L_{\tau(u)}(T_u\tau.\xi)
\end{eqnarray*}

\noindent (and likewise $\theta^L$ changes into $\theta^R$). Since in
addition to that $\tau$ is a group automorphism and an isometry for $(.\,|\,
.)$, the expression (\ref{step}) becomes~:

\begin{eqnarray*}
(\ref{step}) & = & \Big(\theta^L_{\tau\big(\mu_1(x_1)\big)}
\big(T_{\mu_1(x_1)}\tau . (T_{x_1}\mu_1 . v_1)\big) \,|\,
\theta^R_{\tau\big(\mu_2(x_2)\big)} \big(T_{\mu_2(x_2)}\tau . (T_{x_2}\mu_2
.w_2)\big) \Big) \\
& = & \Big(T_1 \tau . \big(\theta^L_{\mu_1(x_1)}(T_{x_1}\mu_1 .v_1)\big) \,|\,
T_1 \tau . \big(\theta^R_{\mu_2(x_2)}(T_{x_2}\mu_2 .w_2)\big) \Big)\\
& = & \big(\theta^L_{\mu_1(x_1)}(T_{x_1}\mu_1 .v_1) \,|\,
\theta^R_{\mu_2(x_2)}(T_{x_2}\mu_2 .w_2) \big)\\
& = & \big( (\mu_1^*\theta^L)_{x_1}(v_1) \,|\, (\mu_2^*\theta^R)_{x_2}(w_2)
\big)
\end{eqnarray*}
\noindent so that we have~:
\begin{eqnarray*}
(\beta^*\w)_x (v,w) & = & (\ref{prefinal}) \\
& = & -(\w_1\oplus\w_2)_x(v,w)\\
& & - \frac{1}{2} \Big( \big(
(\mu_1^*\theta^L)_{x_1}(v_1) \,|\, (\mu_2^*\theta^R)_{x_2}(w_2)\big) - \big(
(\mu_1^*\theta^L)_{x_1}(w_1) \,|\, (\mu_2^*\theta^R)_{x_2}(v_2)\big) \Big)\\
& = & -(\w_1\oplus\w_2)_x(v,w) - (\mu_1^*\theta^L \wedge \mu_2^*\theta^R)_x
(v,w) \\
& = & -\w_x(v,w)
\end{eqnarray*}
\noindent which completes the proof of lemma \ref{anti-inv}.
\end{proof}

\begin{rk} As a last comment on theorem \ref{lag_moduli}, we would like to
say that even if we drop the assumption on the conjugacy classes,
the set of equivalence clases of Lagrangian representations is still a
Lagrangian submanifold of $\Mod$. Indeed, it is always contained in
$Fix(\hat{\beta})$, and therefore it is isotropic, and its dimension is half
the dimension of $\Mod$ (see \cite{FW}). With our hypothesis on the
$\calC_j$, the upshot is that we are able to show that this Lagrangian
submanifold is exactly $Fix(\hat{\beta})$.
\end{rk}
The main tool to obtain theorem \ref{lag_moduli} was theorem \ref{lag_locus}, which is
very general. It may for instance help to find Lagrangian submanifolds in
the moduli space of polygons in $S^3$, which also admits a quasi-Hamiltonian
description (see \cite{Treloar}). In fact, in \cite{Treloar}, the symplectic
structure of the moduli space of polygons with fixed sidelengths in
$S^3 \simeq SU(2)$ is obtained by reduction from the quasi-Hamiltonian space
$\pconj$ where $\calC_j$ is a conjugacy class in $SU(2)$, so that our
involution $\beta$ can be defined in this context. In analogy with results
in \cite{Foth2}, the fixed-point set of this involution should consist of
polygons in $S^3$ which are contained in the equatorial $S^2 \subset S^3$ (I
would like to thank Philip Foth for suggesting this to me).

\subsection{The case of an arbitrary compact connected Lie group}

To conclude, we wish to explain, using the description of $\widehat{U(n)}$
given in section \ref{laginv}, how to make the notion of
Lagrangian representation make sense when the compact connected Lie group $U$ at hand
is not necessarily the unitary group $U(n)$. We suppose that such a group
$U$ is endowed with an involution $\tau$ leaving a \emph{maximal} torus of $U$
pointwise fixed (for instance the Cartan involution
defining its non-compact dual), and we define an action of $\Z / 2\Z =
\{1,\s_0\}$ on $U$ by $\s_0.u=\tau(u)$. We then consider the semi-direct
product $U \rtimes \Z/2\Z$ for this action. Recall that if $U=U(n)$ and
$\tau(u)=\overline{u}$ then $U(n) \rtimes \Z/2\Z = \widehat{U(n)} = U(n)
\sqcup U(n)\s_{L_0}$. Under this identification,
$\s_0.u=\tau(u)=\overline{u}=\s_{L_0}u\s_{L_0}$ and the Lagrangian involutions are the
elements $\s_L=\s_{\pphi(L_0)}=\pphi\s_{L_0}\pphi^{-1}=
(\pphi\s_{L_0}\pphi^{-1}\s_{L_0})\s_{L_0} =
(\pphi\overline{\pphi}^{-1})\s_{L_0} = (\pphi\pphi^t)\s_{L_0} \leftrightarrow
(\pphi\pphi^t,\s_{L_0}) \in U(n)\rtimes\Z/2\Z$. Observe that the element
$\pphi\pphi^t$ does not depend on the choice of $\pphi\in U(n)$ such that
$L=\pphi(L_0)$, as was shown in proposition \ref{angles}. Thus, we see again
that the elements of order $2$ that we are interested in are in one-to-one
correspondence with the symmetric elements of $U(n)$. In the general case, the elements of order $2$ that we are
interested in are the elements $(w,\s_0)\in U \rtimes \Z/2\Z$ where $w\in U$
satisfies $\tau(w)=w^{-1}$. The product of two such elements is then of the
form $(w_1,\s_0).(w_2,\s_0)=(w_1(\s_0.w_2),\s_0^2) = (w_1\tau(w_2),1) \in U
\subset U\rtimes \Z/2\Z$ (observe that when $w_2=w_1$, we indeed obtain $1$
because $\tau(w_1)=w_1^{-1}$). One can then say that a $U$-representation
$(u_1, \,\ldots\, ,u_l)$ of
$\pi=\piS$ is \emph{decomposable} (or \emph{Lagrangian}) if there exist
$w_1, \,{\ldots}\, , w_l \in U$ such that $\tau(w_j)=w_j^{-1}$ for all $j$ and
$u_1=(w_1,\s_0).(w_2,\s_0), u_2=(w_2,\s_0).(w_3,\s_0), \,{\ldots}\, ,
u_l=(w_l,\s_0).(w_1,\s_0)$. Observe that we then have indeed $u_1{\ldots}
u_l=1$, for $u_1{\ldots} u_l =
(w_1\tau(w_2),1).(w_2\tau(w_3),1){\ldots}(w_l\tau(w_1),1) =
(w_1\tau(w_2)w_2\tau(w_3){\ldots}w_l\tau(w_1),1) = 1$ since
$\tau(w_j)=w_j^{-1}$. A representation will be called
$\s_0$\emph{-decomposable} if it is decomposable with $w_1=Id$. Then, theorem
\ref{s0-Lag} and corollary \ref{charac}, along with theorem \ref{lag_moduli} are still
true in this setting (the condition on the eigenvalues of some $\calC_j$ to
be pairwise distinct is to be replaced by the condition that the centralizer
$Z_u$ of any $u \in \calC_j$ is a maximal torus of $U$, and therefore
conjugate to a maximal torus fixed pointwise by $\taum$). All one has to do is then
define $\beta$ as in (\ref{def_beta}) in subsection
\ref{real_problems} (that is, replace $u^t$ by $\tau(u^{-1})$ in the definition of
$\beta$ given in subsection \ref{beta})~: the $\s_0$-decomposable
representations are exactly the elements of the fixed-point set of $\beta$,
a given representation $u$ is decomposable if and only if $\beta(u)$ is equivalent to
$u$, and the set of equivalence classes of decomposable representations is a
Lagrangian submanifold of $\Hom_{\mathcal{C}}(\pi,U)/U$, obtained as the fixed-point set
of an antisymplectic involution $\hat{\beta}$. We refer to \cite{mythesis} for
further details in that direction.

\bibliographystyle{alpha}

\begin{thebibliography}{GHJW97}

\bibitem[AB83]{AB}
M.F. Atiyah and R.~Bott.
\newblock The {Y}ang-{M}ills equations over {R}iemann surfaces.
\newblock {\em Philos. Trans. Roy. Soc. London Ser. A}, 308(1505):523--615,
  1983.

\bibitem[AKSM02]{AKSM}
A.~Alekseev, Y.~Kosmann-Schwarzbach, and E.~Meinrenken.
\newblock Quasi-{P}oisson manifolds.
\newblock {\em Canad. J. Math.}, 54(1):3--29, 2002.

\bibitem[AM95]{AM}
A.~Alekseev and A.~Malkin.
\newblock Symplectic structure of the moduli space of flat connections on a
  {R}iemann surface.
\newblock {\em Comm. Math. Phys.}, 169(1):99--119, 1995.

\bibitem[AMM98]{AMM}
A.~Alekseev, A.~Malkin, and E.~Meinrenken.
\newblock {L}ie group valued moment maps.
\newblock {\em J. of Differential Geom.}, 48(3):445--495, 1998.

\bibitem[AMW01]{AMW}
A.~Alekseev, E.~Meinrenken, and C.~Woodward.
\newblock Linearization of {P}oisson actions and singular values of matrix
  products.
\newblock {\em Ann. Inst. Fourier (Grenoble)}, 51(6):1691--1717, 2001.

\bibitem[Ati82]{Ati}
M.F. Atiyah.
\newblock Convexity and commuting {H}amiltonians.
\newblock {\em Bull. London Math. Soc.}, 14(1):1--15, 1982.

\bibitem[AW98]{AW}
S.~Agnihotri and C.~Woodward.
\newblock Eigenvalues of products of unitary matrices and quantum {S}chubert
  calculus.
\newblock {\em Math. Res. Lett.}, 5(6):817--836, 1998.

\bibitem[Bel01]{Be}
P.~Belkale.
\newblock Local systems on $\mathbb{P}^1 \backslash {S}$ for ${S}$ a finite
  set.
\newblock {\em Compositio Math.}, 129(1):67--86, 2001.

\bibitem[Bis98]{Bi1}
I.~Biswas.
\newblock A criterion for the existence of a parabolic stable bundle of rank
  two over the projective line.
\newblock {\em Internat. J. of Math.}, 9(5):523--533, 1998.

\bibitem[Bis99]{Bi2}
I.~Biswas.
\newblock On the existence of unitary flat connections over the punctured
  sphere with given local monodromy around the punctures.
\newblock {\em Asian J. Math}, 3(2):333--344, 1999.

\bibitem[Dui83]{Du}
J.J. Duistermaat.
\newblock Convexity and tightness for restrictions of {H}amiltonian functions
  to fixed point sets of an antisymplectic involution.
\newblock {\em Trans. Amer. Math. Soc.}, 275(1):417--429, 1983.

\bibitem[EL05]{Lu-Ev}
S.~Evens and J.H. Lu.
\newblock Thompson's conjecture for real semi-simple {L}ie groups.
\newblock In {\em The breadth of symplectic and Poisson geometry}, number 232
  in Prog. Math., pages 121--137. Birkh\"auser, 2005.

\bibitem[Fal01]{F}
E.~Falbel.
\newblock Finite groups generated by involutions on {L}agrangian planes in
  $\mathbb{C}^2$.
\newblock {\em Canad. Math. Bull.}, 44(4):408--419, 2001.

\bibitem[FH05]{Foth2}
P.~Foth and Y.~Hu.
\newblock Toric degenerations of weight varieties and applications.
\newblock In {\em Travaux math\'ematiques. Fasc. XVI}, Trav. Math., XVI, pages 87--105. Univ. Luxembourg, Luxembourg, 2005.

\bibitem[FMS04]{FMS}
E.~Falbel, J.P. Marco, and F.~Schaffhauser.
\newblock Classifying triples of {L}agrangians in a {H}ermitian vector space.
\newblock {\em Topology Appl.}, 144(1-3):1--27, 2004.

\bibitem[Ful98]{Fu}
W.~Fulton.
\newblock Eigenvalues of sums of {H}ermitian matrices (after {A}. {K}lyachko).
\newblock {\em Ast\'erisque}, No. 252(Exp. No. 845):255--269, 1998.
\newblock S\'eminaire Bourbaki.

\bibitem[FW06]{FW}
E.~Falbel and R.~Wentworth.
\newblock Eigenvalues of products of unitary matrices and {L}agrangian
  involutions.
\newblock {\em Topology}, 45(1):55--99, 2006.

\bibitem[Gal97]{gal}
A.J. Galitzer.
\newblock {\em On the moduli space of closed polygonal linkages on the
  2-sphere}.
\newblock PhD thesis, University of Maryland, 1997.

\bibitem[GHJW97]{GHJW}
K.~Guruprasad, J.~Huebschmann, L.~Jeffrey, and A.~Weinstein.
\newblock Group systems, groupoids, and moduli spaces of parabolic pundles.
\newblock {\em Duke Math. J.}, 89(2):377--412, 1997.

\bibitem[Gol84]{G}
W.M. Goldman.
\newblock The symplectic nature of fundamental groups of surfaces.
\newblock {\em Adv. in Math.}, 54(2):200--225, 1984.

\bibitem[Hel01]{He}
S.~Helgason.
\newblock {\em Differential Geometry, Lie Groups and Symmetric Spaces}.
\newblock Number~34 in Graduate Series in Mathematics. AMS, 2001.

\bibitem[Jef94]{J}
L.~Jeffrey.
\newblock Extended moduli spaces of flat connections on {R}iemann surfaces.
\newblock {\em Math. Ann.}, 298(4):667--692, 1994.

\bibitem[JW92]{JW}
L.~Jeffrey and J.~Weitsman.
\newblock {B}ohr-{S}ommerfeld orbits in the moduli space of flat connections
  and the {V}erlinde dimension formula.
\newblock {\em Comm. Math. Phys.}, 150(3):593--630, 1992.

\bibitem[KM99]{KM}
M.~Kapovich and J.~Millson.
\newblock On the moduli space of a spherical polygonal linkage.
\newblock {\em Canad. Math. Bull.}, 42:307--320, 1999.

\bibitem[Knu00]{K}
A.~Knutson.
\newblock The symplectic and algebraic geometry of {Horn}'s problem.
\newblock {\em Linear Algebra Appl.}, 319(1-3):61--81, 2000.

\bibitem[Loo69]{Lo}
O.~Loos.
\newblock {\em Symmetric spaces, II : Compact spaces and classification}.
\newblock W.A. Benjamin, Inc., 1969.

\bibitem[LR91]{Lu-Ra}
J.H. Lu and T.~Ratiu.
\newblock On the nonlinear convexity theorem of {K}ostant.
\newblock {\em J. Amer. Math. Soc.}, 4(2):349--363, 1991.

\bibitem[LS91]{LS}
E.~Lerman and R.~Sjamaar.
\newblock Stratified symplectic spaces and reduction.
\newblock {\em Ann. of Math.}, 134(2):375--422, 1991.

\bibitem[Mor01]{Morita}
S.~Morita.
\newblock {\em Geometry of Differential Forms}.
\newblock Number 201 in Translations of Mathematical Monographs. AMS, 2001.

\bibitem[MW99]{MW}
E.~Meinrenken and C.~Woodward.
\newblock Cobordism for {H}amiltonian loop group actions and flat connections
  on the punctured two-sphere.
\newblock {\em Math. Z.}, 231(1):133--168, 1999.

\bibitem[Nic91]{Nicas}
A.J. Nicas.
\newblock Classifying pairs of {L}agrangians in a {H}ermitian vector space.
\newblock {\em Topology Appl.}, 42(1):71--81, 1991.

\bibitem[OS00]{OSS}
L.~O'Shea and R.~Sjamaar.
\newblock Moment maps and {R}iemannian symmetric pairs.
\newblock {\em Math. Ann.}, 317(3):415--457, 2000.

\bibitem[Sch]{S}
F.~Schaffhauser.
\newblock Un th\'eor\`eme de convexit\'e r\'eel pour les applications moment \`a valeurs dans un groupe de Lie.
\newblock Submitted. Preprint available at http://arxiv.org/abs/math.SG/0609517.

\bibitem[Sch05]{mythesis}
F.~Schaffhauser.
\newblock {\em Decomposable representations and Lagrangian submanifolds of
  moduli spaces associated to surface groups}.
\newblock PhD thesis, Universit\'e Pierre et Marie Curie, 2005. http://www.institut.math.jussieu.fr/theses/2005/schaffhauser/

\bibitem[Tre02]{Treloar}
T.~Treloar.
\newblock The symplectic geometry of polygons in the 3-sphere.
\newblock {\em Canad. J. Math.}, 54(1):30--54, 2002.

\end{thebibliography}

\noindent\footnotesize{\emph{Institut de Math\'ematiques\\
Universit{\'e} Pierre et Marie Curie-Paris 6\\
4, place Jussieu\\
F-75252 Paris Cedex 05\\
email~:} \texttt{florent@math.jussieu.fr}}

\end{document}